\documentclass[11pt]{amsart}

\usepackage[T1]{fontenc}
\usepackage[utf8]{inputenc}
\usepackage{amssymb,mathtools,mathrsfs,booktabs,array}
\usepackage{tikz}
\usepackage[shortalphabetic,initials,nobysame]{amsrefs}
\usepackage[colorlinks=true,linkcolor=red,urlcolor=magenta,citecolor=black]{hyperref}
\hypersetup{
 pdftitle={Calogero-Moser Polynomial Approximations and Dwork Congruences for Complete Flags},
 pdfauthor={Peter Koroteev},
 pdfsubject={Complete-flag vertex functions, Dwork congruences, and tCM operators},
 pdfkeywords={Dwork congruences, complete flags, Calogero-Moser, quasimap vertex functions}
}
\newtheorem{theorem}{Theorem}[section]
\newtheorem{lemma}[theorem]{Lemma}
\newtheorem{proposition}[theorem]{Proposition}
\newtheorem{corollary}[theorem]{Corollary}
\theoremstyle{definition}
\newtheorem{definition}[theorem]{Definition}
\newtheorem{example}[theorem]{Example} 
\theoremstyle{remark}
\newtheorem{remark}[theorem]{Remark}

\numberwithin{equation}{section}
\newcommand{\Z}{\mathbb Z}
\newcommand{\Q}{\mathbb Q}
\newcommand{\C}{\mathbb C}

\newcommand{\F}{\mathbb F}
\DeclareMathOperator{\CT}{CT}
\DeclareMathOperator{\Fl}{Fl}

\newcommand{\eps}{\epsilon}

\newcommand{\sgn}{\operatorname{sgn}}

\renewcommand{\eprint}[1]{\href{https://arxiv.org/abs/#1}{arXiv:#1}}
\BibSpec{article}{
    +{}  {\PrintAuthors}                {author}
    +{,} { \textit}                     {title}
    +{.} { }                            {part}
    +{:} { \textit}                     {subtitle}
    +{,} { \PrintContributions}         {contribution}
    +{.} { \PrintPartials}              {partial}
    +{,} { }                            {journal}
    +{}  { \textbf}                     {volume}
    +{}  { \PrintDatePV}                {date}
    +{,} { \issuetext}                  {number}
    +{,} { \eprintpages}                {pages}
    +{,} { }                            {status}
    +{,} { \eprint}                     {eprint}
    +{}  { \parenthesize}               {language}
    +{}  { \PrintTranslation}           {translation}
    +{;} { \PrintReprint}               {reprint}
    +{.} { }                            {note}
    +{.} {}                             {transition}
    +{}  {\SentenceSpace \PrintReviews} {review}
}
\title[Calogero--Moser Approximations for Complete Flags]{Calogero--Moser Polynomial Approximations\texorpdfstring{\\}{ }and Dwork Congruences for Complete Flags}
\author{Peter Koroteev}
\address{Beijing Institute for Mathematical Sciences and Applications,
Beijing, China}
\email{peter.koroteev@gmail.com}
\date{November 1, 2026}

\begin{document}
\begin{abstract}
We construct $p$-adic polynomial approximations to the zero-framing
cohomological quasimap vertex of cotangent bundles to complete flag varieties from the
trigonometric Calogero--Moser system. Negative integer couplings
give regular staircase Jack polynomials with nonnegative integer
coefficients. We compare them with complete-flag master
polynomials extending the construction of Smirnov--Varchenko.
Both families converge coefficientwise to the vertex and satisfy
the full commuting equations modulo $p^s$. For $n\geq3$, we
prove sharp degree-dependent congruences and compute leading
residues at every $p$-power scale. Their discrepancy has
coefficient norm $p^{-1}$ and an exact weighted Gauss norm giving
uniform decay on smaller closed polydiscs. For the master
family, multivariable Dwork congruences reconstruct the vertex
by an infinite product. The results hold in every rank
for $0<r/d\leq1/2$ and odd prime $p$.
\end{abstract}
\maketitle
\setcounter{tocdepth}{1}
\tableofcontents

\section{Introduction}
\label{sec:introduction}
Enumerative geometry of Nakajima quiver varieties relates quantum
cohomology and quantum $K$-theory to integrable systems
\cites{Nak94,MO19,KPSZ21,K23}. For the cotangent bundle to a
complete flag variety,
\[
 X_n=T^*\Fl(\C^n),\qquad n\geq2,
\]
the $K$-theoretic vertex, with its scalar character factor, is an
eigenfunction of the trigonometric
Ruijsenaars--Schneider (tRS) Hamiltonians in the K\"ahler
parameters \cite{K21,Din21}. Its cohomological limit solves the
trigonometric Calogero--Moser (tCM) equations, which also describe
the quantum differential module of the Springer resolution
\cite{BMO11}.
The recent work \cite{K26} establishes compatibility of the
complete-flag tRS eigenvalue problem with Frobenius action at
roots of unity: partial shifts cancel in Newton--Adams
combinations of the Hamiltonians, leaving the operators obtained
by Frobenius action. Here we use the cohomological tCM spectral
problem to construct polynomial approximations to the vertex
and study their arithmetic properties.

A $K$-theoretic precedent is provided by \cite{K21}
(Proposition~2.11), where a discrete specialization of framing
characters truncates a normalized complete-flag vertex to a
Macdonald polynomial. These polynomials occur in the
correspondence between $A$-type quiver varieties and ADHM moduli
spaces. In the cohomological setting, polynomial solutions arise
by specializing the coupling to negative integers. Choosing these
integers to approach a fixed rational coupling $p$-adically
produces a family of exact tCM eigenfunctions at varying
couplings and polynomial solutions of the fixed system modulo
$p^s$.
After multiplication by a monomial, these terminating
cohomological vertices are symmetric Jack polynomials
\cite{Mac95,NS13}, the differential counterparts of the
symmetric Macdonald polynomials in \cite{K21}.
We use these polynomial specializations to construct arithmetic
approximations of the quasimap vertex. The spectral problem
identifies the terminating vertices with Jack polynomials;
our main arithmetic result compares them with the
master-polynomial approximations.

A different arithmetic construction was introduced by Smirnov
and Varchenko \cite{SV23}. They express a Grassmannian vertex
as a period integral
of a master function, a product of rational powers of differences
of auxiliary variables and K\"ahler parameters. Such period
representations are familiar from mirror symmetry
\cite{Hori03}. The construction here uses the $3d$ mirror
description of the vertex.
To obtain arithmetic approximations, they replace the rational
exponents by nonnegative integers approaching them $p$-adically.
The integrand becomes a master polynomial. Extracting a prescribed
coefficient in its auxiliary variables replaces the period
integral and gives a polynomial in the K\"ahler parameters.
These polynomials satisfy Dwork congruences and reconstruct
the vertex by an infinite product.
We extend this construction to complete flags and prove a sharp
prime-power comparison with the regular negative-parameter
staircase Jack family, uniformly in rank. The two families have
the same coefficientwise limit and satisfy the same full
commuting equations modulo $p^s$. Their difference has a
degree-dependent valuation bound and explicit nonzero leading
residues at every $p$-power scale. These formulas describe the
arithmetic of the Jack coefficients in the monomial symmetric
basis and distinguish the two canonical approximation procedures.

\subsection{Integrability and Polynomial Truncation}
Write $Q=(Q_1,\ldots,Q_{n-1})$ for the K\"ahler parameters
(Novikov variables). The vertex is a generating series of
equivariant quasimap invariants of $X_n$; its coefficient in
multidegree $\boldsymbol d$ records the corresponding curve class.
Quasimaps are maps from a parametrized
rational curve with finitely many base points allowed.
We take the vertex with trivial insertion at the identity
fixed point and normalize the source-rotation weight to one.
The cotangent weight is denoted by $\kappa$. After forming the
complete localization sum in each degree, we specialize the
additive framing weights to zero. We prove that this
specialization is regular and write
\[
 V_n(Q;\kappa),\qquad V_n(0;\kappa)=1,
\]
for the resulting cohomological vertex. Each coefficient is a
polynomial in $\kappa$ over $\Q$.
The tCM equation constructs this series recursively by formal
perturbation theory around the free Hamiltonian. At any fixed
K\"ahler degree, only finitely many perturbative terms contribute.
Specializing the coupling to a negative integer then gives the
terminating polynomials used below; the expansion is given in
Section~\ref{cm:formal-perturbation}.

Let $Q_i=z_i/z_{i+1}$, $D_i=z_i\partial_{z_i}$, and set
\begin{equation}\label{in:calogero}
 H_{n,\kappa}=\frac12\sum_iD_i^2
 -\kappa(\kappa-1)\sum_{i<j}\frac{z_iz_j}{(z_i-z_j)^2},
 \qquad
 J_{n,\kappa}=\prod_{i<j}(1-z_i/z_j)^\kappa.
\end{equation}
All expressions are expanded in the chamber $Q=0$.
The differential tRS limit gives $n$ commuting operators
$\mathsf I_{n,k}(\kappa)$ of order $k$, with principal symbols
$e_k(\xi_1,\ldots,\xi_n)$, where $e_k$ is the elementary
symmetric polynomial of degree $k$.
The first is total momentum; on functions of the ratios $Q_i$,
the second is $-H_{n,\kappa}$. In vertex
gauge we put
\[
 \mathsf B_k^{(0)}(\kappa)
 =J_{n,\kappa}^{-1}\mathsf I_{n,k}(\kappa)J_{n,\kappa}.
\]
The vertex-gauge equation
$\mathcal L_{n,\kappa}V=0$, where
$\mathcal L_{n,\kappa}=J_{n,\kappa}^{-1}H_{n,\kappa}J_{n,\kappa}$,
has a unique normalized formal solution in characteristic zero.
We shall identify it with both the
geometric vertex and the triangular chamber period.
We write $P_\lambda^{(\alpha)}$ for the monic Jack polynomial
with parameter $\alpha=1/\kappa$, and set
$G_n(z)=\prod_{i=1}^n z_i^{i-1}$.

\begin{theorem}[Theorem~\ref{con:vertex} and
Propositions~\ref{pt:polynomiality}, \ref{pt:jack-realization},
\ref{pt:positive-support}, \ref{pt:joint-comparison}]
\label{in:polynomial-theorem}
For every $m\in\Z_{\geq0}$,
\[
 F_m^{(n)}(Q):=V_n(Q;-m)\in\Z[Q_1,\ldots,Q_{n-1}],
 \qquad F_m^{(n)}(0)=1.
\]
Their coefficients are nonnegative integers. For $m\geq1$,
$G_n^mF_m^{(n)}=P_{m\delta}^{(-1/m)}$, where
$\delta=(n-1,\ldots,1,0)$.
These polynomials satisfy the exact equations
\[
 \mathsf B_k^{(0)}(-m)F_m^{(n)}=0,\qquad 1\leq k\leq n.
\]
Fix $\omega=r/d$ in lowest terms, $0<\omega\leq\tfrac12$,
and an odd prime $p$. Put
\[
 P_s=p^s,\qquad A_s=\omega(P_s-1),\qquad
 \widetilde T_s^{(n)}=F_{A_s}^{(n)}.
\]
Then $\widetilde T_s^{(n)}$ converges coefficientwise
$p$-adically to $V_n(Q;\omega)$ and
\[
 \mathsf B_k^{(0)}(\omega)\widetilde T_s^{(n)}
 \equiv0\pmod{p^s},\qquad 1\leq k\leq n.
\]
\end{theorem}

Equivalently, the vertex is the coefficientwise $p$-adic limit
\[
 V_n(Q;\omega)
 =\lim_{s\to\infty}G_n(z)^{-A_s}
    P_{A_s\delta}^{(-1/A_s)}(z).
\]
Both the partition and the Jack polynomial parameter vary with $s$; the
normalization makes their limit a series in the fixed K\"ahler
variables. The exact spectral equations at coupling $-A_s$
and the constant coefficient one uniquely determine this
polynomial family in characteristic zero.
The specialization differs from the
framing specialization in \cite{K21}; the common principle is
to obtain polynomial eigenfunctions from an enumerative series.
For $n=2$ this is the termination of
${}_2F_1(-m,-m;1;Q)$. In higher rank polynomiality follows from
the triangular constant term and bounded external flows.

\subsection{Comparison of the Two Families}
The master-polynomial family is constructed from a triangular
array with $i$ auxiliary variables in row $i$, for
$1\leq i<n$, and the external variables $z_1,\ldots,z_n$ in
the final row. The master function multiplies a factor
$x^{\omega-1}$ for each auxiliary variable $x$, a factor
$(x-y)^{2\omega}$ for each pair in the same row, and a factor
$(x-y)^{-\omega}$ for each pair in consecutive rows. A normalized
integral over interlacing contours represents the vertex.

Following \cite{SV23}, we replace these exponents by nonnegative
integers depending on $p^s$ and approaching the original exponents
$p$-adically. The resulting master polynomial has integral
coefficients. Extracting a prescribed coefficient in its auxiliary
variables and removing an external monomial gives
$T_s^{(n)}$, normalized by $T_s^{(n)}(0)=1$ and $T_0^{(n)}=1$.
Definition~\ref{ar:approximants} specifies the exponents and
coefficient extraction.
The two constructions approximate the same vertex, but their
polynomials need not coincide. Define their difference by
\[
 \Delta_s^{(n)}=\widetilde T_s^{(n)}-T_s^{(n)}.
\]
For a rational number $c\ne0$, let $v_p(c)$ be the exponent of
$p$ in $c$ and set $|c|_p=p^{-v_p(c)}$, with $|0|_p=0$
and $v_p(0)=+\infty$.
The field $\Q_p$ is the completion of $\Q$ for this norm, and
$\Z_p=\{c\in\Q_p:|c|_p\leq1\}$ is its ring of integers.
On $\Z_p[[Q]]$ the coefficient norm is
$\|\sum c_{\boldsymbol e}Q^{\boldsymbol e}\|_p
=\sup_{\boldsymbol e}|c_{\boldsymbol e}|_p$.
A congruence modulo $p^s$ means divisibility of every
coefficient. Coefficientwise convergence is weaker than
convergence in this norm.

\begin{theorem}[Theorems~\ref{pt:coefficientwise-bound} and
\ref{pt:sharpness}]\label{in:comparison-theorem}
Under the arithmetic hypotheses above, both families converge
coefficientwise to $V_n(Q;\omega)$ and satisfy its full commuting
tCM equations modulo $p^s$. Let $n\geq3$, $s\geq1$, and
put $a_0=\omega(p-1)$.
For every $1\leq i\leq n-2$ and $0\leq\ell<s$,
\begin{equation}\label{in:sharp-residue}
 [Q_i^{p^\ell}Q_{i+1}^{p^\ell}]\Delta_s^{(n)}
 \equiv p^{s-\ell}a_0^2\pmod{p^{s-\ell+1}}.
\end{equation}
In particular,
\begin{equation}\label{in:comparison-norm}
 \|\widetilde T_s^{(n)}-T_s^{(n)}\|_p=p^{-1}
\end{equation}
at every level, although every fixed coefficient of the
difference tends to zero.
\end{theorem}

The full degree-dependent estimate is given in
Theorem~\ref{pt:coefficientwise-bound}. The displayed residues
attain this estimate at every $p$-power scale and occur on every
consecutive three-particle subflag.

The comparison follows from the constant-term formulas for the
two families. For a fixed K\"ahler degree, the constant-term
conditions bound the binomial indices that can contribute.
Frobenius divisibility of the resulting binomial coefficients
yields the degree-dependent valuation bound. Coefficients
supported on three consecutive particles attain this bound,
establishing its sharpness in every rank. Since the proof
requires no division by orbit sizes, it also applies to
admissible primes $p\leq n$.

Multiplication by $G_n^{A_s}$ gives a congruence in the ring of
symmetric polynomials. Set $m=A_s$,
$\mathcal T_s=G_n^mT_s^{(n)}$, and $\delta=(n-1,\ldots,1,0)$.
Then
\[
 P_{m\delta}^{(-1/m)}-\mathcal T_s
 \in p\Z[z_1,\ldots,z_n]^{S_n}.
\]
Both polynomials are monic with the same permutahedral Newton
polytope. Formula~\eqref{in:sharp-residue} computes the
coefficient difference obtained by transferring $p^\ell$ units
from one part of $m\delta$ to the part two positions later.
The resulting partition is a three-row squeeze of the staircase;
see Corollary~\ref{pt:symmetric-comparison}. These formulas give
explicit coefficient congruences comparing the regular
negative-parameter staircase Jack family with its
master-polynomial approximants.

The weighted Gauss norm gives a quantitative comparison between
coefficientwise and norm convergence. For $0<\varrho\leq1$,
let $\|\cdot\|_{p,\varrho}$ denote the Gauss norm in which the
coefficient of $Q^{\boldsymbol e}$ is weighted by
$\varrho^{|\boldsymbol e|}$. We prove
\[
 \|\Delta_s^{(n)}\|_{p,\varrho}
 =\max_{0\leq\ell<s}p^{-(s-\ell)}\varrho^{2p^\ell}.
\]
For every fixed $\varrho<1$ this is $O_\varrho(p^{-s})$;
at radius one it remains $p^{-1}$
(Corollary~\ref{pt:weighted-gauss}). Thus the discrepancy
converges uniformly on every smaller closed polydisc, while
its coefficient norm stays fixed. For either family, the
distance to the vertex itself is
\begin{equation}\label{in:distance-to-vertex}
 \|T_s^{(n)}-V_n\|_p
 =\|\widetilde T_s^{(n)}-V_n\|_p=1.
\end{equation}
This follows from a simple-root restriction and Lucas
congruences; see Corollary~\ref{pt:distance-to-limit}.
The ratios $T_s^{(n)}(Q)/T_{s-1}^{(n)}(Q^p)$, called
Frobenius quotients, converge in coefficient norm.
Here $Q^p=(Q_1^p,\ldots,Q_{n-1}^p)$.

\subsection{Dwork Congruences and the Commuting Equations}
Smirnov and Varchenko \cite{SV23} proved Grassmannian Dwork
congruences and anticipated extensions to partial flags.
Their scalar ghost argument already uses alternation and a
Vandermonde expansion. We verify the row-degree bounds and
chamber normalization for the complete-flag dimension vector
$1,\ldots,n-1$. Kartik and Smirnov \cite{KaS25} construct a
related $q$-deformation in the Grassmannian case.

\begin{theorem}[Theorems~\ref{ar:limit}, \ref{ar:dwork},
\ref{ar:product}, and~\ref{tc:full-congruence}]
\label{in:arithmetic-theorem}\label{in:product-theorem}
\label{in:energy-theorem}
Under the arithmetic hypotheses above, for every $n\geq2$,
the polynomials $T_s^{(n)}\in\Z[Q]$
have constant coefficient one and converge coefficientwise to
$V_n(Q;\omega)$. For $s\geq1$,
\begin{equation}\label{in:dwork}
 T_{s+1}^{(n)}(Q)T_{s-1}^{(n)}(Q^p)
 \equiv T_s^{(n)}(Q)T_s^{(n)}(Q^p)\pmod{p^s}.
\end{equation}
For every $m\geq1$,
\begin{equation}\label{in:product}
 V_n(Q;\omega)\equiv
 \prod_{j\geq0}
 \frac{T_m^{(n)}(Q^{p^j})}{T_{m-1}^{(n)}(Q^{p^{j+1}})}
 \pmod{p^m}.
\end{equation}
The product is taken in the $(Q_1,\ldots,Q_{n-1})$-adic topology.
Moreover, for $1\leq k\leq n$,
\begin{equation}\label{in:joint-congruences}
 \mathsf B_k^{(0)}(\omega)V_n=0,\qquad
 \mathsf B_k^{(0)}(\omega)T_s^{(n)}\equiv0\pmod{p^s}.
\end{equation}
All statements hold without a restriction on rank relative to $p$.
\end{theorem}

The commuting system has two roles in the construction. It
constructs the exact-coupling polynomials in
Theorem~\ref{in:polynomial-theorem} and identifies the geometric
vertex with the triangular period. Its integral coefficients
also permit the transfer of differential equations to arithmetic
approximants. The latter step is formulated as a general principle.
Write $T_s=T_s^{(n)}$ and $V=V_n(Q;\omega)$; the ratio
$T_s/V$ is a formal series because $V(0)=1$.
Under the Dwork congruences,
\[
 \theta^\alpha(T_s/V)\in p^s\Z_p[[Q]],\qquad \alpha\ne0.
\]
Here $\theta_i=Q_i\partial_{Q_i}$ is the Euler derivative in
the $i$th K\"ahler variable, and
$\theta^\alpha=\prod_i\theta_i^{\alpha_i}$ for a multi-index
$\alpha=(\alpha_1,\ldots,\alpha_{n-1})$ of nonnegative integers.
Equivalently, every positive-order Euler derivative of $T_s/V$
vanishes modulo $p^s$.
It follows that every integral finite-order differential
annihilator of $V$ annihilates $T_s$ modulo $p^s$
(Lemma~\ref{tc:dwork-transfer}). The full commuting congruences
therefore follow from this principle and the integral tCM
operators.
The Dwork congruences proved here concern $T_s^{(n)}$; no
Dwork congruence for the exact-coupling Jack family is asserted.
The comparison theorem provides information not determined by
these equations. The joint equations, Dwork congruences, and
normalization do not uniquely determine the approximants, since
they admit Frobenius modifications; see
Remark~\ref{pt:frobenius-ambiguity}.

At $\omega=1/2$, the symmetry of the physical coupling
$\kappa(\kappa-1)$ gives joint congruences modulo $p^{2s}$
for the exact-coupling eigenfunctions in their
natural Schr\"odinger gauge; see \eqref{pt:half-coupling-full}.

\subsection{Two- and Three-Particle Examples}
For $n=2$, the vertex is the Gauss hypergeometric series
$V_2(Q;\omega)={}_2F_1(\omega,\omega;1;Q)$. The two
polynomial families coincide:
\[
 T_s^{(2)}(Q)=\widetilde T_s^{(2)}(Q)
 =\sum_{d=0}^{A_s}\binom{A_s}{d}^{\!2}Q^d.
\]
For $\omega=1/2$ and $p=3$, one has $T_1^{(2)}=1+Q$,
and the product formula reads
$V_2(Q;1/2)\equiv\prod_{j\geq0}(1+Q^{3^j})\pmod3$.
This is the projective-line example of \cite{SV23}.
The two constructions first differ for $n=3$.

For $n=3$, $\omega=1/2$, and $p=3$, the first level is
\begin{align*}
 T_1^{(3)}&=(1+Q_1+Q_1Q_2)(1+Q_2+Q_1Q_2),\\
 \widetilde T_1^{(3)}&=T_1^{(3)}+3Q_1Q_2.
\end{align*}
At the next level, the families differ modulo $9$:
\[
 [Q_1^3Q_2^3]\Delta_2^{(3)}=136704\equiv3\pmod9.
\]
This coefficient occurs on every consecutive three-particle
subflag in higher rank.

Formula~\eqref{in:product} determines the normalized quasimap
coefficients modulo $p^m$ from two finite polynomials.
Coefficients of total degree at most $D$ require only $K$ factors
with $p^K>D$. For example, at $p=3$ and precision $9$,
\[
 V_3(Q;1/2)\equiv
 \frac{T_2^{(3)}(Q)T_2^{(3)}(Q^3)}{T_1^{(3)}(Q^3)}
 \pmod{9\Z_3[[Q]]+(Q_1,Q_2)^9}.
\]
The coefficients considered here are zero-framing specializations
of equivariant quasimap invariants. The variables remain formal;
no numerical $p$-adic convergence of $V_n$ is required.

\subsection{Structure of the Paper}
Section~\ref{ar:section-master} defines the triangular period and
the master-polynomial family. Section~\ref{cg:section} develops
the commuting tCM equations, explains their perturbative
solution, and constructs the exact-coupling polynomials.
Section~\ref{pt:section} proves the sharp
comparison. The Dwork and transfer arguments follow in
Section~\ref{ar:section-ghosts}, and
Section~\ref{fl:section} records the explicit calculations.
Appendix~\ref{dg:section}
derives the localization formula and the cohomological normalization.
Appendix~\ref{tc:section} gives the kernel identities in the
integral chamber normalization used here.

\subsection{Acknowledgements}
I thank Alexander Varchenko for discussions of Dwork
congruences and the initial complete-flag calculations. I also
thank Andrey Smirnov for useful discussions of vertex functions and
their arithmetic properties. I acknowledge support from the Beijing Natural Science Foundation
(BJNSF) under grant IS26008.

I used ChatGPT 6.1 iteratively in order to check the calculations for complete flag varieties and higher Dwork congruences, polish the proofs, and refine the statements.

\section{Triangular Master Polynomials and Chamber Periods}
\label{ar:section-master}
We use the triangular quiver of dimension vector $1,\ldots,n-1$
to construct a chamber period and the master-polynomial family.
Section~\ref{cg:section} identifies the period with the zero-framing
quasimap vertex and constructs the exact-coupling polynomial
approximations from the tCM equations.
The localization convention is recorded in Appendix~\ref{dg:section}.

The period construction follows \cite{SV23} (\S3). Its mirror
symmetry interpretation replaces an enumerative generating
function by an integral of an explicit multivalued density.
Period representations are standard in mirror symmetry
\cite{Hori03}; here the relevant construction is supplied by
$3d$ mirror symmetry and verified by the normalized tCM equation.
The quiver determines the factors of the density, and a choice
of interlacing contours determines its expansion in the
K\"ahler variables.

Let $n\geq2$ and $N=\binom n2$. Write the internal variables as
$x=(x_{i,a})_{1\leq a\leq i\leq n-1}$ and set
$x_{n,a}=z_a$, $1\leq a\leq n$, for the final row. Define
\begin{align}
 M(x)&=\prod_{i=1}^{n-1}\prod_{a=1}^i x_{i,a},
 &\Delta(x)&=\prod_{i=1}^{n-1}\prod_{1\leq a<b\leq i}
                   (x_{i,b}-x_{i,a}),\label{ar:M-Delta}\\
 E(x,z)&=\prod_{i=1}^{n-1}\prod_{a=1}^i\prod_{b=1}^{i+1}
                   (x_{i,a}-x_{i+1,b}),
 &G_n(z)&=\prod_{a=1}^n z_a^{a-1},
 \qquad c_n=\binom{n+1}{3}.\label{ar:E-G}
\end{align}
We use $\mathbf d_{i,a}=a$ and the all-ones vector $\mathbf1$ of
length $N$. The coefficient of
$\prod_{i,a}x_{i,a}^{\nu_{i,a}}$ in $F$ is denoted by
$[x^{\boldsymbol\nu}]F$; for a polynomial this coefficient
vanishes if a component of $\boldsymbol\nu$ is negative. The
chamber coordinates are
\begin{equation}\label{ar:Q}
 Q_j=z_j/z_{j+1},\qquad 1\leq j\leq n-1.
\end{equation}

Choose distinct positive radii for the internal variables and the
$z_a$, satisfying
\begin{equation}\label{ar:interlacing}
 |x_{i+1,a}|<|x_{i,a}|<|x_{i+1,a+1}|
 \qquad(1\leq a\leq i\leq n-1).
\end{equation}
These inequalities admit such choices and imply
$|z_1|<\cdots<|z_n|$. For variables in adjacent rows we denote
the ratio of smaller to larger radius by $\rho_{i,a,b}$:
\begin{equation}\label{ar:edge-ratio}
 \rho_{i,a,b}=
 \begin{cases}
 x_{i+1,b}/x_{i,a},& b\leq a,\\
 x_{i,a}/x_{i+1,b},& b>a.
 \end{cases}
\end{equation}
The adjacent-row ratios and the same-row ratios
$x_{i,a}/x_{i,b}$, $a<b$, all have absolute value less than one
on this torus.

\begin{definition}\label{ar:period}
For $\omega\in\Q$, define the normalized chamber period by
\begin{equation}\label{ar:period-CT}
 \mathcal V_n(Q;\omega)=\CT_x
 \left[
 \prod_{i=1}^{n-1}\prod_{a<b}
       (1-x_{i,a}/x_{i,b})^{2\omega}
 \prod_{i=1}^{n-1}\prod_{a=1}^i\prod_{b=1}^{i+1}
       (1-\rho_{i,a,b})^{-\omega}
 \right].
\end{equation}
We expand each factor by the binomial theorem in the indicated
ratio. The notation $\CT_xF=[x^{\mathbf0}]F$ extracts the
coefficient of $\prod_{i,a}x_{i,a}^0$ from this Laurent expansion
and retains the external variables $z_a$.
\end{definition}

For $n=2$, the chamber is $|z_1|<|x|<|z_2|$ and
\[
 \CT_x(1-z_1/x)^\lambda(1-x/z_2)^\mu
 =\sum_{k\geq0}\binom\lambda k\binom\mu k Q_1^k.
\]
The equality of the two indices is conservation of flow at $x$.
The same interpretation controls the triangular period.

\begin{lemma}\label{ar:flow}
For arbitrary rational exponents in the binomial factors of
\eqref{ar:period-CT}, the prescribed constant term belongs to
$1+(Q_1,\ldots,Q_{n-1})\Q[[Q_1,\ldots,Q_{n-1}]]$.
Each coefficient receives finitely many contributions from the
binomial expansion.
\end{lemma}

\begin{proof}
Assign a graph vertex to every internal and external variable,
and an edge to every binomial factor. We orient each edge from smaller to larger radius.
The binomial indices are nonnegative integer flows on the edges;
the exponent at a vertex is outgoing minus incoming flow.
Constant-term extraction is therefore conservation of flow at
the internal vertices. The ordering of radii makes the graph
acyclic.

Let $\mu_a$ be the exponent of $z_a$ in a surviving monomial.
The sum of all vertex exponents vanishes, giving $\sum_a\mu_a=0$
and
\[
 z_1^{\mu_1}\cdots z_n^{\mu_n}
 =\prod_{j=1}^{n-1}Q_j^{e_j},
 \qquad e_j=\sum_{a=1}^j\mu_a.
\]
Consider a cut between consecutive radii in the gap from $|z_j|$
to $|z_{j+1}|$. Flow conservation at internal vertices identifies
the sum of its crossing flows with $e_j$. Since all crossing edges
have the same orientation, $e_j\geq0$.

Every internal radius is between $|z_1|$ and $|z_n|$.
Consequently each edge crosses a cut in one of the external gaps,
and its flow is bounded by the corresponding $e_j$, hence by
$\max_j e_j$. The graph has finitely many edges, so a fixed
multidegree receives finitely many contributions. In degree zero
all flows vanish, and the product of the binomial constant terms
is one.
\end{proof}

The argument extends the chamber expansion of \cite{SV23}
(Proposition~3.1 and the proof of Theorem~4.4) to complete flags.
It also separates the integrality of the period from the
additional hypotheses needed for the polynomial approximants.

\begin{proposition}\label{ar:general-integrality}
Every coefficient of $\mathcal V_n(Q;\omega)$ is a polynomial
in $\omega$ over $\Q$. For every prime $p$ and rational
$\omega\in\Z_p$,
\[
 \mathcal V_n(Q;\omega)\in\Z_p[[Q_1,\ldots,Q_{n-1}]].
\]
In particular, if $\omega=r/d$ in lowest terms with $d>0$, then
\[
 \mathcal V_n(Q;r/d)\in\Z[1/d][[Q_1,\ldots,Q_{n-1}]].
\]
\end{proposition}

\begin{proof}
Lemma~\ref{ar:flow} expresses each coefficient as a finite sum
of products of binomial polynomials in $2\omega$ and $-\omega$.
This proves polynomial dependence on $\omega$. If
$\omega\in\Z_p$, both upper arguments belong to $\Z_p$ and
all their binomial coefficients are $p$-integral. Indeed, for
fixed $k$, the polynomial $\binom{\alpha}{k}$ is continuous on
$\Z_p$ and takes integer values on the dense subset
$\Z_{\geq0}$. The finite sum is therefore integral. Applying
this to every prime not dividing $d$ proves the last assertion.
\end{proof}

The constant-term formula is equivalent to a contour integral
of the multivalued density
\begin{equation}\label{ar:master-density}
 \mathcal M_\omega(x,z)=
 M(x)^{\omega-1}\Delta(x)^{2\omega}E(x,z)^{-\omega}.
\end{equation}
Factor each difference by the variable of larger radius and use
the branches specified by the binomial expansions. The $c_n$
negative adjacent-row factors contribute a constant phase.
With $(-t)^{-\omega}=e^{-\pi\mathrm i\omega}t^{-\omega}$, the
normalized differential form becomes
\begin{equation}\label{ar:normalized-density}
 e^{\pi\mathrm i\omega c_n}G_n(z)^\omega
 \mathcal M_\omega(x,z)\bigwedge_{i,a}dx_{i,a}
 =\left[\text{the product in \eqref{ar:period-CT}}\right]
   \bigwedge_{i,a}\frac{dx_{i,a}}{x_{i,a}}.
\end{equation}
The wedge products have the same ordering. In contour notation,
\begin{equation}\label{ar:period-contour}
 \mathcal V_n(Q;\omega)
 =\frac{e^{\pi\mathrm i\omega c_n}G_n(z)^\omega}
        {(2\pi\mathrm i)^N}
   \oint_{\Gamma(z)}
   M(x)^{\omega-1}\Delta(x)^{2\omega}E(x,z)^{-\omega}
   \bigwedge_{i,a}dx_{i,a}.
\end{equation}
Here $\Gamma(z)$ is the product of positively oriented circles
with the radii in \eqref{ar:interlacing}; branches are fixed by
\eqref{ar:normalized-density}. The normalized integrand is
single-valued on this torus. Integrating its chamber expansion
extracts precisely the internal constant term.
In particular, \eqref{ar:master-density} is a density for the
measure $\bigwedge dx_{i,a}$; its expression as a constant term
already accounts for the logarithmic measure.

For the master-polynomial construction we impose
\begin{equation}\label{ar:prime-assumptions}
 \omega=r/d,\qquad (r,d)=1,\qquad 0<\omega\leq\tfrac12,
 \qquad p\text{ odd}.
\end{equation}
For $s\geq0$ put
\begin{equation}\label{ar:exponents}
\begin{gathered}
 P_s=p^s,\qquad A_s=\omega(P_s-1),\qquad
 \beta_s=(1-\omega)(P_s-1),\\
 \gamma_s=1+(1-2\omega)(P_s-1)=P_s-2A_s.
\end{gathered}
\end{equation}
The exponents are nonnegative integers, with $\gamma_s$ odd.
The master polynomial is
\begin{equation}\label{ar:Phi}
 \Phi_s(x,z)=M(x)^{\beta_s}\Delta(x)^{\gamma_s}E(x,z)^{A_s}.
\end{equation}
At $s=0$, we have $\beta_0=A_0=0$ and $\gamma_0=1$.
Thus $\Phi_0=\Delta$, the product of the rowwise Vandermonde
factors. The group $\mathfrak S=\prod_{i=1}^{n-1}S_i$ acts by
permutations within each internal row: $S_i$ is the symmetric group
on $x_{i,1},\ldots,x_{i,i}$, and all external variables $z_a$ are
fixed. The odd exponent $\gamma_s$ makes $\Phi_s$ alternating
under this action.

\begin{definition}\label{ar:approximants}
Set
\begin{equation}\label{ar:T-def}
 U_s(z)=[x^{P_s\mathbf d-\mathbf1}]\Phi_s(x,z),\qquad
 T_s^{(n)}(Q)=(-1)^{A_sc_n}G_n(z)^{-A_s}U_s(z).
\end{equation}
The coefficient extracted in \eqref{ar:T-def} is that of
$\prod_{i,a}x_{i,a}^{ap^s-1}$ in $\Phi_s$; it is a polynomial
in the external variables $z_a$. The sign and monomial factors
normalize its constant coefficient. The following proposition
shows that the resulting expression is an integral polynomial
in $Q$ with constant coefficient one.
\end{definition}

\begin{proposition}\label{ar:normalization}
One has $T_s^{(n)}\in\Z[Q_1,\ldots,Q_{n-1}]$ and
$T_s^{(n)}(0)=1$, including $T_0^{(n)}=1$. More explicitly,
\begin{equation}\label{ar:T-CT}
 T_s^{(n)}=\CT_x
 \left[
 \prod_{i=1}^{n-1}\prod_{a<b}
       (1-x_{i,a}/x_{i,b})^{\gamma_s}
 \prod_{i=1}^{n-1}\prod_{a=1}^i\prod_{b=1}^{i+1}
       (1-\rho_{i,a,b})^{A_s}
 \right].
\end{equation}
\end{proposition}

\begin{proof}
We factor the differences in \eqref{ar:Phi} by their larger-radius
variables. The variable $x_{i,a}$ has $a-1$ smaller variables in
its own row and $2a-1$ in the adjacent rows, interpreting the
absent first row as empty. Its monomial exponent is therefore
\[
 \beta_s+(a-1)\gamma_s+(2a-1)A_s=aP_s-1.
\]
The exponent of $z_a$ is $(a-1)A_s$. The sign is determined by
the adjacent-row factors whose first variable has smaller radius.
Their number is
\[
 \sum_{i=1}^{n-1}\sum_{a=1}^i(i+1-a)
 =\sum_{i=1}^{n-1}\frac{i(i+1)}2=c_n.
\]
Hence the monomial and its sign are
\begin{equation}\label{ar:leading-monomial}
 (-1)^{A_sc_n}G_n(z)^{A_s}x^{P_s\mathbf d-\mathbf1}.
\end{equation}
Dividing by this monomial converts the selected coefficient into
the constant term \eqref{ar:T-CT}. Its nonnegative integer powers
give a finite expansion over $\Z$. Lemma~\ref{ar:flow} then shows
that the result is a polynomial in $Q$ and has constant coefficient
one.
\end{proof}

\begin{lemma}\label{ar:master-support}
For $P=p^s$, a variable on row $i$ has degree at most
\[
 \deg_{x_{i,a}}\Phi_s\leq iP+A_s-1<(i+1)P-1.
\]
A coefficient $[x^{P\boldsymbol u-\mathbf1}]\Phi_s$ can be
nonzero only when $\boldsymbol u=\tau(\mathbf d)$ for a rowwise
permutation $\tau\in\mathfrak S$. Its value is then
$\sgn(\tau)U_s(z)$.
\end{lemma}

\begin{proof}
A variable on row $i$ occurs in $i-1$ same-row differences and
$2i$ adjacent-row differences. Thus its degree is bounded by
$\beta_s+(i-1)\gamma_s+2iA_s=iP+A_s-1$.
Since $A_s<P$, an exponent $Pu_{i,a}-1$ must have
$1\leq u_{i,a}\leq i$. Alternation makes the coefficient zero
if two entries in a row agree. The remaining entries are a
permutation of $1,\ldots,i$, with the asserted sign.
\end{proof}

\begin{example}\label{ar:rank-one}
For $n=2$, the coefficient extraction is straightforward:
\[
 T_s^{(2)}(Q;\omega)=\sum_{k=0}^{A_s}\binom{A_s}{k}^{\!2}Q^k.
\]
After removing the monomial factor in $x_{1,1}$, we extract the
coefficient of $x_{1,1}^{A_s}$ in
$(x_{1,1}-z_1)^{A_s}(x_{1,1}-z_2)^{A_s}$ and use the normalization
of \eqref{ar:T-def}. The limiting period is
\[
 \mathcal V_2={}_2F_1(\omega,\omega;1;Q).
\]
This is the projective-line example of \cite{SV23}
(Example~2.2).
\end{example}

Lemma~\ref{ar:flow} also controls the passage to the $p$-adic
limit. In each fixed K\"ahler degree the binomial indices remain
bounded, so we can take the limit in their finite sum.

\begin{theorem}\label{ar:limit}
Under \eqref{ar:prime-assumptions},
\[
 \lim_{s\to\infty}[Q^{\boldsymbol e}]T_s^{(n)}(Q)
 =[Q^{\boldsymbol e}]\mathcal V_n(Q;\omega)
 \quad\text{in }\Q_p
\]
for every $\boldsymbol e\in\Z_{\geq0}^{n-1}$. In particular,
$\mathcal V_n\in\Z_p[[Q_1,\ldots,Q_{n-1}]]$.
\end{theorem}

\begin{proof}
By Lemma~\ref{ar:flow}, all binomial indices in multidegree
$\boldsymbol e$ are at most $\max_j e_j$. The coefficient in
\eqref{ar:T-CT} is therefore a finite sum of products of binomial
polynomials in $A_s$ and $\gamma_s$, with an indexing set
independent of $s$. An index exceeding a nonnegative integer
upper argument contributes zero. The limits $A_s\to-\omega$
and $\gamma_s\to2\omega$ in $\Q_p$ give \eqref{ar:period-CT} by
continuity of these rational polynomials. The approximating
coefficients belong to $\Z$, and the closedness of $\Z_p$ in
$\Q_p$ proves that their limits are integral.
\end{proof}

The polynomials $T_s^{(n)}$ therefore approximate the period
$p$-adically in every fixed K\"ahler degree.
Theorem~\ref{con:vertex} identifies this period
with the quasimap partition function, so the same construction
proves its $p$-integrality. Increasing $s$ improves each fixed-degree
coefficient; the polynomials $T_s^{(n)}$ need not be Taylor
truncations of the period. Theorem~\ref{ar:limit} is the
complete-flag counterpart of \cite{SV23} (Theorem~4.4).

\section{Integrability and Exact Polynomial Approximations}
\label{cg:section}
The tCM system gives a second construction of polynomial
approximations to the vertex. Formal perturbation theory around
the free Hamiltonian constructs the normalized Schr\"odinger
wavefunction. Passing to vertex gauge gives a series whose
coefficients depend polynomially on the coupling. At negative
integer coupling the vertex series terminates to an integral
polynomial.
Choosing these integers to approach $\omega$ in $\Q_p$ produces
exact eigenfunctions at varying couplings and congruences for the
fixed system at $\omega$.

As in \cite{K21} (Proposition~2.11), we obtain polynomial
eigenfunctions by specializing parameters of an enumerative
series. Here we specialize the cohomological coupling at zero
additive framing. We prove termination and its arithmetic
consequences directly from the normalized chamber solution.

\subsection{The Vertex Gauge and Formal Uniqueness}
Put
\begin{equation}\label{cg:variables}
 Q_i=\frac{z_i}{z_{i+1}},\qquad
 q_{ij}=\frac{z_i}{z_j}=Q_i\cdots Q_{j-1},\qquad
 u_{ij}=\frac{q_{ij}}{1-q_{ij}}\quad(i<j).
\end{equation}
For $D_i=z_i\partial_{z_i}$ and
$\theta_i=Q_i\partial_{Q_i}$, with $\theta_0=\theta_n=0$,
\begin{equation}\label{cg:relative-derivatives}
 D_i=\theta_i-\theta_{i-1},\qquad
 H_0:=\frac12\sum_iD_i^2
 =\sum_{i=1}^{n-1}\theta_i^2
   -\sum_{i=1}^{n-2}\theta_i\theta_{i+1}
\end{equation}
on ratio functions. In particular, $\sum_iD_i=0$ on these
functions. For a coupling $\kappa$, use the Hamiltonian normalization
\cite{C71,S71,M75}
\begin{equation}\label{cg:schrodinger}
\begin{gathered}
 H_{n,\kappa}=H_0-g_\kappa U_n,\qquad
 g_\kappa=\kappa(\kappa-1),\qquad
 U_n=\sum_{i<j}\frac{q_{ij}}{(1-q_{ij})^2},\\
 J_{n,\kappa}=\prod_{i<j}(1-q_{ij})^\kappa.
\end{gathered}
\end{equation}
Every factor is expanded at $Q=0$ with constant coefficient one.
We abbreviate $J_{n,\omega}=J_n$.

\begin{proposition}\label{cg:gauge}
The vertex-gauge Hamiltonian is
\begin{equation}\label{cg:normalized-operator}
 \mathcal L_{n,\kappa}:=J_{n,\kappa}^{-1}H_{n,\kappa}J_{n,\kappa}
 =H_0-\kappa\sum_{i<j}u_{ij}(D_i-D_j)
      -\kappa^2\sum_{i<j}(j-i)u_{ij}.
\end{equation}
\end{proposition}

\begin{proof}
Let $\alpha_{ij}=e_i-e_j$. The logarithmic gradient of the gauge
is $-\kappa\sum_{i<j}\alpha_{ij}u_{ij}$. Its divergence gives
$-\kappa\sum_{i<j}q_{ij}/(1-q_{ij})^2$ in the conjugated kinetic
term. The remaining multiplication term is
\begin{equation}\label{cg:zero-order-intermediate}
 \kappa^2\left[
 \frac12\left|\sum_{i<j}\alpha_{ij}u_{ij}\right|^2
 -\sum_{i<j}(u_{ij}+u_{ij}^2)\right].
\end{equation}
For $i<j<k$, the relation $q_{ik}=q_{ij}q_{jk}$ gives
\[
 u_{ij}u_{ik}-u_{ij}u_{jk}+u_{ik}u_{jk}=-u_{ik}.
\]
Grouping cross terms by triples therefore yields
\[
 \frac12\left|\sum_{i<j}\alpha_{ij}u_{ij}\right|^2
 =\sum_{i<j}u_{ij}^2-\sum_{i<j}(j-i-1)u_{ij}.
\]
Substitution proves the formula.
\end{proof}

\begin{proposition}\label{cg:formal-uniqueness}
Over a characteristic-zero field containing $\kappa$, there is a
unique series $F_{n,\kappa}\in1+(Q_1,\ldots,Q_{n-1})K[[Q]]$
satisfying $\mathcal L_{n,\kappa}F_{n,\kappa}=0$.
Its coefficients belong to $\Q[\kappa]$ before specialization.
\end{proposition}

\begin{proof}
For $\boldsymbol d\geq0$, put $d_0=d_n=0$. The diagonal
coefficient in degree $\boldsymbol d$ is
\begin{equation}\label{cg:indicial-form}
 \mathcal E(\boldsymbol d)
 =\sum_{i=1}^{n-1}d_i^2-\sum_{i=1}^{n-2}d_id_{i+1}
 =\frac12\sum_{i=1}^n(d_i-d_{i-1})^2.
\end{equation}
It is a positive integer for $\boldsymbol d\ne0$.
All other operator coefficients have positive chamber degree.
Starting from the constant coefficient one, the equation thus
determines every coefficient by induction on total degree.
The only divisions are by $\mathcal E(\boldsymbol d)$, so each
coefficient is polynomial in $\kappa$ over $\Q$.
\end{proof}

For $n=2$, with $Q=z_1/z_2$ and $\theta=Q\partial_Q$,
\begin{equation}\label{cm:two-particle-vertex}
 (1-Q)\mathcal L_{2,\kappa}
 =\theta^2-Q(\theta+\kappa)^2.
\end{equation}
Consequently
\begin{equation}\label{cm:hypergeometric-wave}
 F_{2,\kappa}(Q)={}_2F_1(\kappa,\kappa;1;Q),\qquad
 H_{2,\kappa}\bigl((1-Q)^\kappa F_{2,\kappa}\bigr)=0.
\end{equation}
This is the two-particle equation of \cite{SV23} (Example~2.2).
The chamber solution in higher rank belongs to the type-$A$
Calogero--Moser hypergeometric framework; see \cite{HR15} for
the commuting operators and their recursive kernels.

\begin{remark}\label{cg:uniqueness-scope}
Characteristic-zero uniqueness does not determine an integral
solution modulo $p^s$. The diagonal
$\mathcal E(\boldsymbol d)$ can vanish modulo $p$, and Euler
derivatives annihilate every series in $Q^p$ modulo $p$.
The polynomial construction and the comparison theorem below
retain information that the reduced differential equations do
not specify.
\end{remark}

\subsection{Formal Perturbation Theory}
\label{cm:formal-perturbation}
The normalized chamber solution can be constructed by perturbing
the free Hamiltonian $H_0$ in the physical coupling
$g_\kappa=\kappa(\kappa-1)$. On the ideal
$\mathfrak m=(Q_1,\ldots,Q_{n-1})\Q[[Q]]$, the operator
$H_0$ is invertible: its eigenvalue on $Q^{\boldsymbol d}$ is
the positive integer $\mathcal E(\boldsymbol d)$ of
\eqref{cg:indicial-form}. Since $U_n$ has no constant term,
the operator
\[
 \mathsf K_n:\Q[[Q]]\longrightarrow\mathfrak m,
 \qquad f\longmapsto H_0^{-1}(U_nf)
\]
raises total K\"ahler degree. The normalized Schr\"odinger
equation $H_{n,\kappa}\psi_\kappa=0$ is therefore equivalent to
\begin{equation}\label{cm:perturbation-series}
 \psi_\kappa
 =1+g_\kappa\mathsf K_n\psi_\kappa
 =\sum_{j\geq0}g_\kappa^j\mathsf K_n^j1.
\end{equation}
Only $j\leq|\boldsymbol d|$ contributes to the coefficient of
$Q^{\boldsymbol d}$. Thus this is an exact identity in the
$\mathfrak m$-adic topology for every characteristic-zero
coupling; it requires no analytic smallness assumption on
$g_\kappa$.

Theorem~\ref{con:vertex} identifies the vertex-gauge solution as
\[
 V_n(Q;\kappa)=J_{n,\kappa}^{-1}\psi_\kappa.
\]
Specializing $\kappa=-m$ recovers $F_m^{(n)}$.
Proposition~\ref{pt:polynomiality} proves that its K\"ahler
series terminates, and the full commuting equations follow from
the same vertex identification. Perturbation theory therefore
provides a recursive construction of the polynomial family;
the constant-term formula proves its integrality and finite
support.

For two particles the first terms are
\[
 \psi_\kappa(Q)
 =1+g_\kappa Q+\frac{g_\kappa(g_\kappa+2)}4Q^2+O(Q^3).
\]
At $\kappa=-1$, one has $g_\kappa=2$ and
$\psi_{-1}=(1+Q)/(1-Q)$; returning to vertex gauge gives
$F_1^{(2)}=1+Q$.
The expansion also shows that $\psi_\kappa=\psi_{1-\kappa}$,
because it depends only on $g_\kappa$. Consequently,
\begin{equation}\label{cm:coupling-reflection}
 V_n(Q;\kappa)=J_{n,1-2\kappa}V_n(Q;1-\kappa).
\end{equation}
This symmetry explains the higher precision at half coupling
in \eqref{pt:half-coupling-full}. The perturbative recursion
is over $\Q$; the constant-term construction supplies the
integral coefficients required for the arithmetic congruences.

\subsection{The Full Commuting System}
The higher equations follow from the differential limit of the
K\"ahler-variable tRS system. A generating determinant
describes the full family: its elementary coefficients are the
commuting Hamiltonians. At zero coupling it reduces to
$\prod_i(u-D_i)$, so its leading differential terms are the
elementary symmetric polynomials in the momenta. The
$\kappa$-dependent shifts below come from centering the tRS
operators before taking their differential limit. Let
$\Delta_+(z)=\prod_{i<j}(z_j-z_i)$ and define
\begin{equation}\label{tc:determinant}
\begin{aligned}
 \mathscr S_n(u;\kappa)
 &=\Delta_+(z)^{-1}
  \det\left[z_i^{j-1}\bigl(u-D_i-\kappa(j-i)\bigr)\right]_{i,j=1}^n
 \\
 &=\sum_{k=0}^n(-1)^ku^{n-k}\mathsf B_k^{(0)}(\kappa),
 \qquad \mathsf B_0^{(0)}=1.
\end{aligned}
\end{equation}
The determinant is ordered with coefficients preceding
derivatives; factors from distinct rows commute. It is the
chamber normalization of the Sekiguchi determinant \cite{NS13}.
Set
\[
 \mathsf I_{n,k}(\kappa)
 =J_{n,\kappa}\mathsf B_k^{(0)}(\kappa)J_{n,\kappa}^{-1}.
\]
We suppress $n$ and $\kappa$ when they are fixed.

\begin{proposition}\label{dg:full-family}
The operators $\mathsf B_k^{(0)}(\kappa)$ commute and have
principal symbols $e_k(\xi_1,\ldots,\xi_n)$. They are
algebraically independent in the full $n$-variable differential
operator algebra. For generic additive framing $a$, the vertex
in the conventions of Appendix~\ref{dg:section} satisfies
\begin{align}\label{dg:geometric-joint-equations}
 \mathsf B_k^{(0)}(\kappa)
    \bigl(z^{-a}V_n(Q;a,\kappa)\bigr)
 &=e_k(-a)z^{-a}V_n(Q;a,\kappa),\\
 \mathsf I_{n,k}(\kappa)\psi_a&=e_k(-a)\psi_a,
 \qquad \psi_a=z^{-a}J_{n,\kappa}V_n(Q;a,\kappa).\notag
\end{align}
These equations specialize regularly to $a=0$.
\end{proposition}

\begin{proof}
Use the centered commuting tRS family
$\mathsf A_{r,\eta}$ from \eqref{dg:normalized-trs}, with
$q=e^\eta$ and $h=e^{\eta\kappa}$. Expansion by a subset of
rows gives
\begin{equation}\label{cm:generating-limit}
 \sum_{r=0}^n(-1)^re^{-\eta ur}\mathsf A_{r,\eta}
 =\Delta_+(z)^{-1}
 \det\left[z_i^{j-1}
 \left(1-e^{\eta(D_i+\kappa(j-i)-u)}\right)\right]_{i,j=1}^n.
\end{equation}
Indeed, the alternant for the chosen rows factors into the
corresponding tRS subset coefficient; the centering removes its
chamber character. Each matrix entry is divisible by $\eta$.
Dividing by $\eta^n$ and taking the constant coefficient in
$\eta$ gives \eqref{tc:determinant}. The left sides commute for
different $u$, so their regular limits commute. The top
differential part is $\prod_i(u-D_i)$, giving the principal
symbols and algebraic independence.

Proposition~\ref{dg:kahler-input} gives the eigenvalue of the left
side on $z^{-a}\widehat V_n^K$ as
\[
 \sum_{r=0}^n(-1)^re^{-\eta ur}e_r(e^{-\eta a})
 =\prod_i\left(1-e^{-\eta(u+a_i)}\right).
\]
Its quotient by $\eta^n$ tends to $\prod_i(u+a_i)$.
The coefficientwise cohomological limit of
Proposition~\ref{cl:generic-limit} proves the first equation.
Conjugation proves the second, and
Proposition~\ref{dg:vertex-regularity} permits $a=0$.
\end{proof}

The first two operators are
\begin{equation}\label{cm:joint-generators}
\begin{gathered}
 \mathsf B_1^{(0)}=\sum_iD_i,\qquad
 \mathsf B_2^{(0)}=\frac12\left(\sum_iD_i\right)^2
                         -\mathcal L_{n,\kappa},\\
 \mathsf I_{n,2}=\frac12\left(\sum_iD_i\right)^2-H_{n,\kappa}.
\end{gathered}
\end{equation}
Thus $H_{n,\kappa}\psi_a=\frac12\sum_i a_i^2\psi_a$.
The reciprocal spectral characters in Appendix~\ref{dg:section}
account for the signs $e_k(-a)$.

For comparison with the type-$A$ hypergeometric framework, put
$\rho_i=(n+1-2i)/2$ and
\[
 \delta_n(z)=\prod_{i<j}\frac{z_j-z_i}{\sqrt{z_iz_j}}
 =z^{-\rho}\prod_{i<j}(1-q_{ij}).
\]
Removing the Weyl-denominator gauge from $\psi_a$ gives
\[
 \delta_n^{-\kappa}\psi_a
 =z^{\kappa\rho-a}V_n(Q;a,\kappa).
\]
This is the normalized Harish-Chandra chamber series with spectral
parameter $-a$ in the chamber used here; see \cite{HR15} for its
relation to Heckman--Opdam hypergeometric functions. The
negative-integral-coupling constructions of \cite{FV09} supply the
Baker--Akhiezer background. We use the zero-spectral specialization
and prove its polynomiality and integral arithmetic directly.

\begin{theorem}\label{con:vertex}
For every $n\geq2$,
\[
 V_n(Q;\kappa)=\mathcal V_n(Q;\kappa)=F_{n,\kappa}(Q)
 \quad\text{in }\Q[\kappa][[Q]].
\]
In particular, every coupling specialization in characteristic
zero is defined, and
\[
 \mathsf B_k^{(0)}(\kappa)V_n(Q;\kappa)=0,\qquad 1\leq k\leq n.
\]
\end{theorem}

\begin{proof}
The zero-framing vertex satisfies the quadratic equation by
Proposition~\ref{dg:vertex-equation}. The adjacent-rank kernel
identity proves the same equation for the master period in
Proposition~\ref{tc:period-equation}. Both have constant
coefficient one, so Proposition~\ref{cg:formal-uniqueness}
identifies them. Their coefficients are polynomials in the
coupling by Proposition~\ref{ar:general-integrality} and
Proposition~\ref{dg:vertex-regularity}, and the joint equations
follow from Proposition~\ref{dg:full-family}.
\end{proof}

\subsection{Integral Dependence on the Coupling}
We need an integral form of the operators to move the exact
equations between congruent couplings.

\begin{lemma}\label{tc:integral-operators}
Put $\mathscr D_n(Q)=\prod_{i<j}(1-q_{ij})$. For $\tau\in S_n$,
let
\[
 h_a(\tau)=\sum_{i=1}^a\tau(i)-\frac{a(a+1)}2,
 \qquad 1\leq a<n.
\]
Then
\begin{equation}\label{tc:permutation-operator}
 \mathscr D_n(Q)\mathsf B_k^{(0)}(\kappa)
 =\sum_{\tau\in S_n}\sgn(\tau)Q^{h(\tau)}
 e_k\bigl(D_i+\kappa(\tau(i)-i):1\leq i\leq n\bigr).
\end{equation}
Consequently the normal-ordered coefficients belong to
$\Z[\kappa][[Q]]$. For $\lambda,\omega\in\Z_p$ with
$\lambda\equiv\omega\pmod{p^s}$,
\[
 \mathsf B_k^{(0)}(\lambda)
 \equiv\mathsf B_k^{(0)}(\omega)\pmod{p^s}.
\]
\end{lemma}

\begin{proof}
Expand \eqref{tc:determinant} by permutations, with derivatives
to the right. Dividing the alternant monomial by its chamber
leading monomial gives
$\prod_i z_i^{\tau(i)-i}=Q^{h(\tau)}$.
Every $h_a(\tau)$ is nonnegative because the smallest sum of
$a$ distinct integers from $1,\ldots,n$ is $a(a+1)/2$.
The remaining Vandermonde factor is $\mathscr D_n$, an integral
unit with constant coefficient one. This proves the formula
and polynomial dependence on $\kappa$. Specialization at
congruent couplings proves the last assertion.
\end{proof}

In Schr\"odinger gauge, the same operators admit an expansion
over matchings. If
$v_{ij}=q_{ij}/(1-q_{ij})^2$, then
\begin{equation}\label{tc:matching-formula}
 \mathsf I_{n,k}(\kappa)
 =\sum_{\substack{\mathfrak m\text{ a matching},\quad
        S\subseteq\{1,\ldots,n\}\setminus V(\mathfrak m)\\
        2|\mathfrak m|+|S|=k}}
 g_\kappa^{|\mathfrak m|}
 \prod_{\{i,j\}\in\mathfrak m}v_{ij}\prod_{a\in S}D_a.
\end{equation}
Here a matching consists of disjoint unordered pairs, and
$V(\mathfrak m)$ is their set of endpoints. Each matching and
subset occurs once. The formula is the matching expression of
\cite{HR15}; its equality with \eqref{tc:determinant} is proved
in Proposition~\ref{tc:matching-identification}. It shows that
the Schr\"odinger operators have integral chamber coefficients
polynomial in the physical coupling $g_\kappa$.

\subsection{Polynomial Solutions at Negative Integer Coupling}
\label{pt:negative-coupling-construction}
We construct a second family of integral polynomial
approximations by specializing the coupling in the normalized
vertex. Equivalently, we specialize the perturbative solution
\eqref{cm:perturbation-series} and return to vertex gauge.
Using the chamber ratios of
Section~\ref{ar:section-master}, put
\begin{equation}\label{pt:two-exponent-period}
 C_n(Q;m,\gamma)=\CT_x\left[
 \prod_{i=1}^{n-1}\prod_{a<b}
       (1-x_{i,a}/x_{i,b})^\gamma
 \prod_{i=1}^{n-1}\prod_{a=1}^i\prod_{b=1}^{i+1}
       (1-\rho_{i,a,b})^m\right].
\end{equation}
Lemma~\ref{ar:flow} makes every fixed-degree coefficient a finite
sum. Theorem~\ref{con:vertex} identifies
\[
 V_n(Q;\kappa)=C_n(Q;-\kappa,2\kappa).
\]
For $m\in\Z_{\geq0}$, define the exact-coupling solution
\begin{equation}\label{pt:negative-integer-polynomial}
 F_m^{(n)}(Q)=V_n(Q;-m)=C_n(Q;m,-2m).
\end{equation}
The normalized tCM equation fixes this specialization uniquely.

\begin{proposition}\label{pt:polynomiality}
For every integer $m\geq0$,
$F_m^{(n)}\in\Z[Q_1,\ldots,Q_{n-1}]$ and $F_m^{(n)}(0)=1$.
For admissible $\omega,p$ and $A_s=\omega(p^s-1)$,
\[
 F_{A_s}^{(n)}\longrightarrow V_n(Q;\omega)
 \quad\text{coefficientwise in }\Q_p.
\]
\end{proposition}

\begin{proof}
For $m\geq1$, expand
\[
 (1-t)^{-2m}=\sum_{k\geq0}\binom{2m+k-1}{k}t^k.
\]
These same-row coefficients are integers, and the adjacent-row
factors have nonnegative integer exponent $m$.
Lemma~\ref{ar:flow} bounds all expansion indices contributing to
a fixed K\"ahler degree, so its coefficient is a finite sum of
integers and has nonnegative degree.

Each external variable belongs to $n-1$ adjacent-row factors,
whose indices are at most $m$. If its exponent is $\mu_a$,
then $|\mu_a|\leq(n-1)m$. The K\"ahler exponents satisfy
$0\leq e_j=\sum_{a\leq j}\mu_a\leq j(n-1)m$.
Hence only finitely many degrees occur. In degree zero every
flow vanishes, giving constant coefficient one. The case $m=0$
is immediate.

At fixed degree, the period coefficient is polynomial in the
coupling over $\Q$, because the flow bounds are independent of
that coupling. Since $-A_s=\omega(1-p^s)\to\omega$ in
$\Q_p$, these polynomials give the asserted limit.
\end{proof}

The two-particle case is
$F_m^{(2)}(Q)=\sum_{j=0}^m\binom mj^2Q^j$.
In higher rank \eqref{pt:negative-integer-polynomial} provides
both a finite polynomial and its exact joint equations. The
arithmetic use of this construction is to let the coupling vary
through the sequence $-A_s$.

The terminating vertices have a symmetric-polynomial realization.
We use finite-variable monic Jack polynomials in the convention
$\alpha=1/\kappa$ of \cite{Mac95,NS13}. We write $m_\mu$
for the sum of distinct monomials with exponent partition $\mu$.

\begin{proposition}\label{pt:jack-realization}
Let $n\geq2$, $m\geq1$, and
$\delta=(n-1,n-2,\ldots,0)$. The Jack polynomial
$P_{m\delta}^{(\alpha)}(z_1,\ldots,z_n)$ is regular at
$\alpha=-1/m$, and
\begin{equation}\label{pt:jack-specialization}
 G_n(z)^m F_m^{(n)}(z_1/z_2,\ldots,z_{n-1}/z_n)
 =P_{m\delta}^{(-1/m)}(z_1,\ldots,z_n).
\end{equation}
\end{proposition}

\begin{proof}
Since $D_iG_n^\kappa=\kappa(i-1)G_n^\kappa$, conjugation of
\eqref{tc:determinant} gives
\begin{equation}\label{pt:jack-sekiguchi}
 \mathscr S_n^{\mathrm{sym}}(u;\kappa)
 :=G_n^{-\kappa}\mathscr S_n(u;\kappa)G_n^\kappa
 =\Delta_+(z)^{-1}
   \det\left[z_i^{j-1}(u-D_i-\kappa(j-1))\right]_{i,j=1}^n.
\end{equation}
This is the Sekiguchi family on symmetric polynomials. Its
eigenvalue on $P_\mu^{(1/\kappa)}$ is
\[
 \prod_{i=1}^n\bigl(u-\mu_i-\kappa(n-i)\bigr);
\]
see \cite{NS13} (equations~(13) and~(16)), after reversing columns
and making the affine change of spectral variable
$u\mapsto u/\kappa-(n-1)$, with an overall factor $\kappa^n$.

We first justify specialization at the negative parameter.
Write $K_1(\kappa),K_2(\kappa)$ for the first two elementary
coefficients of \eqref{pt:jack-sekiguchi}. The quadratic
operator $K_1(\kappa)^2/2-K_2(\kappa)$ is triangular on
monomial symmetric polynomials, with diagonal
\[
 q_\mu(\kappa)
 =\frac12\sum_{i=1}^n\bigl(\mu_i+\kappa(n-i)\bigr)^2.
\]
Put $\lambda=m\delta$. On the finite dominance ideal spanned
by $m_\mu$ with $\mu\leq\lambda$ and $\ell(\mu)\leq n$,
all operator coefficients are polynomial in $\kappa$.
At $\kappa=-m$,
\[
 q_\lambda(-m)=0,\qquad
 q_\mu(-m)=\frac12\sum_{i=1}^n(\mu_i-\lambda_i)^2>0
 \quad(\mu\ne\lambda).
\]
Thus every denominator $q_\lambda-q_\mu$ in the monic
triangular eigenvector recursion is nonzero at $-m$.
For generic $\kappa$ that eigenvector is
$P_\lambda^{(1/\kappa)}$, so its coefficients are regular
there. The full commuting eigenvalue equations specialize as
well. Their common generating eigenvalue is $u^n$.

It remains to check the chamber normalization. Set
$\gamma=(0,m,\ldots,(n-1)m)$, so $G_n^m=z^\gamma$.
Monicity gives coefficient one at $z^\gamma$, since
$\gamma$ is a permutation of $\lambda$.
For every monomial $z^\beta$ occurring in the Jack polynomial,
the decreasing rearrangement $\mu=\beta^+$ satisfies
$\mu\leq\lambda$. Consequently, for $1\leq j<n$,
\[
 \sum_{i=1}^j\beta_i
 \geq |\lambda|-\sum_{i=1}^{n-j}\mu_i
 \geq |\lambda|-\sum_{i=1}^{n-j}\lambda_i
 =\frac{mj(j-1)}2
 =\sum_{i=1}^j\gamma_i.
\]
Since $|\beta|=|\gamma|$, the quotient
$z^{\beta-\gamma}$ is $Q^{\boldsymbol e}$ with
$e_j=\sum_{i=1}^j(\beta_i-\gamma_i)\geq0$.
Hence $G_n^{-m}P_\lambda^{(-1/m)}$ is a polynomial in $Q$
with constant coefficient one.
By \eqref{pt:jack-sekiguchi}, it satisfies the joint
vertex-gauge equations at $\kappa=-m$, and therefore the
quadratic equation. Proposition~\ref{cg:formal-uniqueness}
identifies it with $F_m^{(n)}$.
\end{proof}

\begin{remark}\label{pt:jack-regularity-scope}
Jack polynomials at negative rational parameters also occur in
the theory of admissible partitions and clustering conditions
\cite{FJMM02,CDL22}. In the convention used here, the coprime
admissibility regime of \cite{FJMM02} is
\[
 \alpha=-\frac{k+1}{r-1},\qquad k\geq1,\quad r\geq2,
 \quad \gcd(k+1,r-1)=1.
\]
The parameter $\alpha=-1/m$ lies outside this regime: it would
require $r-1=m(k+1)$, contrary to coprimality. Thus that
admissibility theorem does not directly give the staircase
regularity used here. Proposition~\ref{pt:jack-realization}
establishes it by the quadratic spectral gap on the finite
dominance ideal. This verifies directly the general noncollision
criterion for regularity; cf.\ \cite{FJMM02} (Lemma~2.4).
The argument concerns the specified finite-variable staircase
polynomial.

The identification itself follows from the classical Sekiguchi
characterization and is consistent with the Macdonald
specialization in \cite{K21} (Proposition~2.11).
The chamber period expresses the coefficients of this family
as sums of nonnegative integers. Its arithmetic comparison with the
master-polynomial approximants is the main result of the
following section.
\end{remark}

The same period presentation proves positivity and determines
the common Newton polytope of the two symmetric polynomial
families.

\begin{proposition}\label{pt:positive-support}
For every $m\geq0$,
\[
 F_m^{(n)}\in\Z_{\geq0}[Q_1,\ldots,Q_{n-1}].
\]
For $m\geq1$, the Jack polynomial in
\eqref{pt:jack-specialization} has nonnegative integer monomial
coefficients. At an arithmetic level, put $m=A_s$ and define
\begin{equation}\label{pt:symmetric-master}
 \mathcal T_s(z)=G_n(z)^mT_s^{(n)}(Q)
 =(-1)^{mc_n}U_s(z).
\end{equation}
Then $\mathcal T_s$ is a symmetric integral homogeneous
polynomial of degree $m\binom n2$. Both $\mathcal T_s$ and
$P_{m\delta}^{(-1/m)}$ have leading monomial symmetric
function $m_{m\delta}$ with coefficient one, and their Newton
polytope is
\[
 \operatorname{conv}\{\tau(m\delta):\tau\in S_n\}.
\]
Consequently, for $1\leq j<n$,
\begin{equation}\label{pt:exact-coordinate-degrees}
 \deg_{Q_j}F_m^{(n)}=m j(n-j),\qquad
 \deg_{Q_j}T_s^{(n)}=m j(n-j).
\end{equation}
\end{proposition}

\begin{proof}
In the constant term for $F_m^{(n)}$, the same-row factors
$(1-t)^{-2m}$ have nonnegative integer coefficients. An
adjacent-row factor contributes $(-1)^{k_e}\binom m{k_e}$.
For a surviving flow, weight the exponent-conservation equation
at a vertex on row $r$ by $r$. Internal conservation and
$\sum_a\mu_a=0$ for the external exponents give
\[
 \sum_{e\text{ joining adjacent rows}}\varepsilon_e k_e=0,
 \qquad \varepsilon_e\in\{1,-1\}.
\]
Same-row edges contribute zero to this weighted sum.
Reducing modulo two shows that $\sum_e k_e$ is even. Thus
every surviving contribution has nonnegative sign. Finiteness
of the sum in each degree proves positivity; $F_0=1$.
Multiplication by $G_n^m$ preserves coefficients, so
Proposition~\ref{pt:jack-realization} gives Jack positivity.

The master polynomial is invariant under permutation of the
external variables. Hence $U_s$, and therefore $\mathcal T_s$,
is symmetric. Since $T_s^{(n)}$ depends only on ratios,
\eqref{pt:symmetric-master} makes its homogeneous degree
$m\binom n2$. Both homogeneous polynomials have coefficient
one at $G_n^m=z^\gamma$, where
$\gamma=(0,m,\ldots,(n-1)m)$. Symmetry gives coefficient one
at every permutation of $m\delta$.

For either polynomial, every occurring exponent vector $\beta$
satisfies
\[
 \sum_{a=1}^j\beta_a\geq\frac{mj(j-1)}2,
 \qquad 1\leq j<n,
\]
because division by $G_n^m$ has nonnegative K\"ahler degrees.
Symmetry gives the same inequality for every subset of $j$
coordinates. Taking complements gives the upper subset-sum
inequalities for the permutahedron of $m\delta$.
Thus the support is contained in this polytope, and every vertex
occurs with coefficient one. This proves the Newton-polytope
statement and dominance triangularity. Finally, maximizing
\[
 e_j=\sum_{a=1}^j\beta_a-\frac{mj(j-1)}2
\]
over the permutation vertices gives $m j(n-j)$, proving
\eqref{pt:exact-coordinate-degrees}.
\end{proof}

\subsection{Arithmetic Specialization and Joint Congruences}
Assume \eqref{ar:prime-assumptions}, set $P=p^s$ and
$m=A_s=\omega(P-1)$, and write
\begin{equation}\label{pt:two-families}
 T_s^{(n)}=C_n(Q;m,P-2m),\qquad
 \widetilde T_s^{(n)}=F_{A_s}^{(n)},\qquad
 \Delta_s^{(n)}=\widetilde T_s^{(n)}-T_s^{(n)}.
\end{equation}
The first equality is the constant-term form of the master
polynomial in Section~\ref{ar:section-master}. Thus the two
families have the same adjacent-row factors, and their same-row
exponents differ by $P$. At $\omega=1/2$, the same-row
exponent of $T_s^{(n)}$ is one. Both constant coefficients are
one, so $\Delta_s^{(n)}(0)=0$.

Put
\[
 \mathcal A_n=\sum_{i<j}u_{ij}(D_i-D_j),\qquad
 \mathcal B_n=\sum_{i<j}(j-i)u_{ij}.
\]
Then $\mathcal L_{n,\kappa}=H_0-\kappa\mathcal A_n
-\kappa^2\mathcal B_n$.

\begin{proposition}\label{pt:energy-comparison}\label{pt:joint-comparison}
For $m\geq0$ and $1\leq k\leq n$,
\begin{equation}\label{pt:exact-joint-polynomial}
 \mathsf B_k^{(0)}(-m)F_m^{(n)}=0.
\end{equation}
In particular,
\begin{equation}\label{pt:shifted-equation}
 \mathcal L_{n,-m}F_m^{(n)}=0,\qquad
 H_{n,-m}(J_{n,-m}F_m^{(n)})=0.
\end{equation}
At the arithmetic level $m=A_s$, one has
\begin{equation}\label{pt:fixed-coupling-residual}
 \mathcal L_{n,\omega}F_m^{(n)}
 =-\omega P\left[\mathcal A_n+(\omega-m)\mathcal B_n\right]F_m^{(n)}
 \in P\Z_p[[Q]],
\end{equation}
and
\begin{equation}\label{pt:joint-congruences}
 \mathsf B_k^{(0)}(\omega)\widetilde T_s^{(n)}\equiv0\pmod{p^s},
 \qquad
 \mathsf B_k^{(0)}(\omega)\Delta_s^{(n)}\equiv0\pmod{p^s}.
\end{equation}
Equivalently, the operators $\mathsf I_{n,k}(\omega)$ annihilate
$J_n\widetilde T_s^{(n)}$ and $J_n\Delta_s^{(n)}$ modulo $p^s$.
\end{proposition}

\begin{proof}
Theorem~\ref{con:vertex} and Proposition~\ref{dg:full-family}
give the exact equations at $\kappa=-m$; specialization is valid
because all vertex coefficients are polynomial in the coupling.
For $\lambda=-m=\omega(1-P)$, Lemma~\ref{tc:integral-operators}
gives
$\mathsf B_k^{(0)}(\lambda)\equiv
\mathsf B_k^{(0)}(\omega)\pmod{p^s}$.
Apply this congruence to the integral polynomial $F_m^{(n)}$.
Subtracting the equations for $T_s^{(n)}$ in
Theorem~\ref{tc:full-congruence} gives the equation for the
difference.

For the explicit quadratic residual, use
\[
 \mathcal L_{n,\omega}-\mathcal L_{n,\lambda}
 =-(\omega-\lambda)
    [\mathcal A_n+(\omega+\lambda)\mathcal B_n]
\]
and $\omega-\lambda=\omega+m=\omega P$.
Its coefficients are integral chamber series. Conjugation by
the integral unit $J_n$ gives the Schr\"odinger congruences.
\end{proof}

In particular,
\begin{equation}\label{pt:modular-kernel}
 \mathcal L_{n,\omega}\Delta_s^{(n)}\equiv0\pmod{p^s},\qquad
 H_{n,\omega}(J_n\Delta_s^{(n)})\equiv0\pmod{p^s}.
\end{equation}
A nonzero difference with zero constant coefficient cannot be an
exact normalized-chamber zero-energy solution in characteristic
zero. The sharp comparison in Section~\ref{pt:section} measures
the additional classes allowed by reduction modulo $p^s$.

The natural Schr\"odinger gauge gives higher precision at
$\omega=1/2$. Set
$\phi_s=J_{n,\lambda}F_m^{(n)}$, with $\lambda=-m$. It is an
integral chamber series and satisfies
\begin{align}\label{pt:physical-residual}
 H_{n,\omega}\phi_s&=(g_\lambda-g_\omega)U_n\phi_s,\\
 g_\lambda-g_\omega&=\omega P(1-2\omega+\omega P).
 \notag
\end{align}
At $\omega=1/2$, one has $g_\lambda-g_{1/2}=P^2/4$.
The matching formula \eqref{tc:matching-formula} is polynomial
in $g_\kappa$ with integral chamber coefficients. Since
$\phi_s$ is an exact joint eigenfunction at $\lambda$, it follows
that
\begin{equation}\label{pt:half-coupling-full}
 \mathsf I_{n,k}(1/2)\phi_s\equiv0\pmod{p^{2s}},
 \qquad 1\leq k\leq n.
\end{equation}
This improvement concerns the natural shifted gauge $\phi_s$;
the comparison of $T_s^{(n)}$ and $\widetilde T_s^{(n)}$ uses the
common gauge $J_n$.

\section{Sharp Comparison of the Polynomial Families}
\label{pt:section}
We now compare the master-polynomial approximations with the exact
negative-coupling solutions of Section~\ref{cg:section}.
Both constant terms use the same factors linking consecutive
rows, raised to the power $m=A_s$. For a pair of variables
within one row, however, the exponent is $p^s-2m$ in $T_s^{(n)}$
and $-2m$ in $\widetilde T_s^{(n)}$.
We first prove that every coefficient of their difference is
divisible by $p$. We then bound the expansion indices that can
contribute to a fixed K\"ahler degree; these bounds give higher
divisibility in low degrees. Finally, a restriction to three
consecutive particles produces explicit nonzero coefficients
and proves that the estimates are sharp.
Throughout this section we assume \eqref{ar:prime-assumptions}
and put $P=p^s$ and $m=A_s=\omega(P-1)$.

\subsection{Frobenius Action and Vandermonde Support}
The Frobenius action on the same-row factors compares the two
families modulo $p$. The resulting constant-term calculation is
controlled by the permutation support of the master polynomial,
so the comparison is valid even when $p\leq n$.

Write
\begin{align}
 \mathscr R_n(x)&=\prod_{i=1}^{n-1}\prod_{a<b}
                       (1-x_{i,a}/x_{i,b}),&
 \mathscr E_n(x,z)&=\prod_{i=1}^{n-1}\prod_{a=1}^i
                   \prod_{b=1}^{i+1}(1-\rho_{i,a,b}),
 \label{pt:row-products}
\end{align}
and put $\mathscr P_s=\mathscr R_n^{P-2m}\mathscr E_n^m$.
The two constant terms are
\begin{equation}\label{pt:comparison-constant-terms}
 T_s^{(n)}=\CT_x\mathscr P_s,\qquad
 \widetilde T_s^{(n)}=\CT_x
                      (\mathscr P_s\mathscr R_n^{-P}).
\end{equation}
The notation $\mathscr R_n(x^P)$ means that each internal variable
is raised to the $P$th power. Its inverse is expanded in the
same smaller-to-larger ratios as the master period.

\begin{lemma}\label{pt:integral-extraction}
For every $n\geq2$ and $s\geq1$, one has the integral identity
\begin{equation}\label{pt:frobenius-extraction}
 \CT_x\left[\mathscr P_s\mathscr R_n(x^P)^{-1}\right]
 =\CT_x\mathscr P_s=T_s^{(n)}.
\end{equation}
\end{lemma}

\begin{proof}
Use the vector $\mathbf d_{i,a}=a$ from
Section~\ref{ar:section-master}. The normalization
\eqref{ar:leading-monomial} reads
\begin{equation}\label{pt:normalized-master-polynomial}
 \mathscr P_s=(-1)^{mc_n}G_n(z)^{-m}
                 x^{-P\mathbf d+\mathbf1}\Phi_s(x,z).
\end{equation}
In particular, the polynomial $\Phi_s$ is alternating in each
internal row. Its degree in a variable on row $i$ is at most
$iP+m-1<(i+1)P-1$. Lemma~\ref{ar:master-support} therefore says that
\[
 [x^{P\boldsymbol u-\mathbf1}]\Phi_s\ne0
 \quad\Longrightarrow\quad
 \boldsymbol u=\tau(\mathbf d),\qquad
 \tau\in\prod_{i=1}^{n-1}S_i.
\]

Every monomial of $\mathscr R_n(x^P)^{-1}$ has the form
$x^{P\boldsymbol\nu}$, where, on each row,
$\boldsymbol\nu$ is a nonnegative integral combination of
$\varepsilon_a-\varepsilon_b$, $a<b$; here
$\varepsilon_a$ denotes the $a$th coordinate vector on that row.
Thus its initial partial sums are nonnegative. A nonzero
contribution to the constant term in
\eqref{pt:frobenius-extraction} requires
\[
 \mathbf d-\boldsymbol\nu=\tau(\mathbf d).
\]
For a row of length $i$, this gives, for $1\leq k<i$,
\[
 0\leq\sum_{a=1}^k\nu_{i,a}
   =\sum_{a=1}^k\bigl(a-\tau_i(a)\bigr)\leq0.
\]
The last inequality holds because $1+\cdots+k$ is the smallest
sum of $k$ distinct elements of $\{1,\ldots,i\}$. Every partial
sum vanishes, hence $\boldsymbol\nu=0$ and
$\tau=\mathrm{id}$. Only the constant monomial of
$\mathscr R_n(x^P)^{-1}$ contributes, proving the identity over
$\Z$. All coefficient extractions are finite by the acyclic-flow
argument of Lemma~\ref{ar:flow}.
\end{proof}

Since $P=p^s$, the Frobenius action gives
$\mathscr R_n(x)^P\equiv\mathscr R_n(x^P)\pmod p$.
Both factors have constant coefficient one in the chamber
expansion, so their inverses obey the same congruence. The lemma
therefore proves divisibility of every discrepancy coefficient by
$p$. We next determine its dependence on the K\"ahler degree.

\subsection{Degree-Dependent Congruences}
For a K\"ahler multidegree
$\boldsymbol e=(e_1,\ldots,e_{n-1})\in\Z_{\geq0}^{n-1}$, set
\begin{equation}\label{pt:degree-size}
 B(\boldsymbol e)=\max_{1\leq j\leq n-2}\min(e_j,e_{j+1}).
\end{equation}
This quantity bounds the flow on every same-row edge in an
appropriate interlacing chamber. For $n=3$ it reduces to
$\min(e_1,e_2)$. We use $v_p(0)=+\infty$.

\begin{theorem}\label{pt:coefficientwise-bound}
Assume \eqref{ar:prime-assumptions}, $n\geq3$, and $s\geq1$.
Then
\begin{equation}\label{pt:global-congruence}
 \Delta_s^{(n)}\in p\Z[Q_1,\ldots,Q_{n-1}].
\end{equation}
If $B(\boldsymbol e)=0$, then
$[Q^{\boldsymbol e}]\Delta_s^{(n)}=0$. If
$B(\boldsymbol e)>0$, then
\begin{equation}\label{pt:degree-precision}
 v_p\bigl([Q^{\boldsymbol e}]\Delta_s^{(n)}\bigr)
 \geq\max\left\{1,\ s-\lfloor\log_pB(\boldsymbol e)\rfloor\right\}.
\end{equation}
In particular, the coefficients with $B(\boldsymbol e)<p$ agree
modulo $p^s$, and every fixed coefficient of
$\Delta_s^{(n)}$ tends to zero $p$-adically.
\end{theorem}

\begin{proof}
Polynomiality follows from Proposition~\ref{pt:polynomiality}
and Proposition~\ref{ar:normalization}. By
\eqref{pt:comparison-constant-terms}, the Frobenius congruence
for the inverse of $\mathscr R_n^P$, and
Lemma~\ref{pt:integral-extraction},
\[
 \widetilde T_s^{(n)}
 \equiv\CT_x[\mathscr P_s\mathscr R_n(x^P)^{-1}]
 =T_s^{(n)}\pmod p.
\]
This proves \eqref{pt:global-congruence} in every rank.

For the degree estimate, write
\begin{equation}\label{pt:difference-kernel}
 \Delta_s^{(n)}=\CT_x\left[
 \mathscr E_n^m\mathscr R_n^{-2m}(1-\mathscr R_n^P)\right].
\end{equation}
Choose the interlacing radii more specifically so that
\[
 |z_{a+n-i-1}|<|x_{i,a}|<|z_{a+n-i}|.
\]
These inequalities are compatible with
\eqref{ar:interlacing}: for example, take
$\log|z_j|=j$ and
$\log|x_{i,a}|=a+n-i-\epsilon(n-i)$, where
$0<\epsilon<1/n$. The binomial ratios in the master period
are unchanged by this refinement of the chamber.

A same-row edge joining $x_{i,a}$ to $x_{i,b}$, $a<b$, crosses
a cut in every external gap with index
$a+n-i-1,\ldots,b+n-i-1$. Flow conservation bounds its total
index by
\[
 \min_{a+n-i-1\leq j\leq b+n-i-1}e_j
 \leq B(\boldsymbol e).
\]
Indeed, the interval contains two adjacent gap indices, and the
sum of all flows across a cut in gap $j$ is $e_j$.

Expand $1-\mathscr R_n^P$ as a sum indexed by the same-row
edges. Each nonconstant term contains at least one binomial
factor $\binom Pk$ with $k>0$. Its index contributes to the total
flow on that edge, so $k\leq B(\boldsymbol e)$. If
$B(\boldsymbol e)=0$, no such term can contribute, proving the
vanishing assertion. If $1\leq B(\boldsymbol e)<P$, use
\begin{equation}\label{pt:binomial-valuation}
 v_p\binom Pk=s-v_p(k),\qquad 1\leq k<P.
\end{equation}
To verify the identity, write
$\binom Pk=(P/k)\binom{P-1}{k-1}$. The second factor is a
$p$-adic unit, since
$(1+X)^{P-1}\equiv\sum_{j=0}^{P-1}(-X)^j\pmod p$.
Thus every contributing term is divisible by
$p^{s-\lfloor\log_pB(\boldsymbol e)\rfloor}$; all remaining
binomial coefficients are integral. If
$B(\boldsymbol e)\geq P$, the global congruence supplies the
bound $1$. Together these prove \eqref{pt:degree-precision}.
At each fixed multidegree the bound tends to infinity with $s$.
\end{proof}

Dividing by $p$ gives an integral solution of the full commuting
congruences at one lower precision. Thus the discrepancy gives
joint tCM solutions in every complete flag.

\begin{corollary}\label{pt:normalized-kernel}
For every $n\geq3$, the polynomial
$E_s^{(n)}=\Delta_s^{(n)}/p$ is integral. For $1\leq k\leq n$,
\[
 \mathsf B_k^{(0)}(\omega)E_s^{(n)}\equiv0\pmod{p^{s-1}},\qquad
 \mathsf I_{n,k}(\omega)(J_{n,\omega}E_s^{(n)})
       \equiv0\pmod{p^{s-1}}.
\]
In particular, $\mathcal L_{n,\omega}E_s^{(n)}\equiv0$ and
$H_{n,\omega}(J_{n,\omega}E_s^{(n)})\equiv0\pmod{p^{s-1}}$.
\end{corollary}

\begin{proof}
Divide \eqref{pt:joint-congruences} by $p$, using
\eqref{pt:global-congruence}, and conjugate by $J_{n,\omega}$.
\end{proof}

\subsection{Consecutive Subflags and Sharpness}
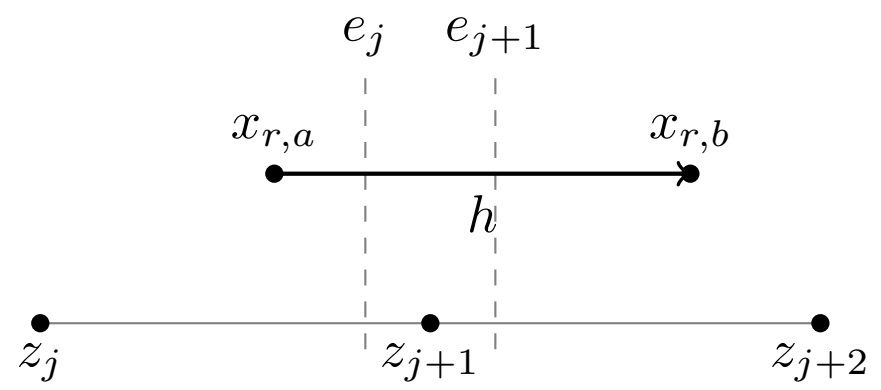
\begin{figure}[t]
\centering
\begin{tikzpicture}[scale=.86,every node/.style={font=\small}]
 \draw[gray] (0,0)--(6,0);
 \foreach \x/\lab in {0/{z_j},3/{z_{j+1}},6/{z_{j+2}}}
  {\fill (\x,0) circle (2pt);\node[below] at (\x,0) {$\lab$};}
 \draw[dashed,gray] (2.5,-.2)--(2.5,1.9);
 \draw[dashed,gray] (3.5,-.2)--(3.5,1.9);
 \node[above] at (2.5,1.9) {$e_j$};
 \node[above] at (3.5,1.9) {$e_{j+1}$};
 \fill (1.8,1.15) circle (2pt);
 \fill (5,1.15) circle (2pt);
 \node[above] at (1.8,1.15) {$x_{r,a}$};
 \node[above] at (5,1.15) {$x_{r,b}$};
 \draw[->,thick] (1.8,1.15)--node[below] {$h$}(5,1.15);
\end{tikzpicture}
\caption{A same-row flow crosses cuts in two adjacent external
gaps. Conservation gives $h\leq\min(e_j,e_{j+1})$.
The refined interlacing chamber in the proof of
Theorem~\ref{pt:coefficientwise-bound} ensures that every
same-row edge crosses such a pair of cuts.}
\label{pt:cut-figure}
\end{figure}

We prove that the coefficient bound is attained at every
$p$-power scale. Restricting to two consecutive K\"ahler
variables reduces the constant term to a three-particle triangle;
its extreme flow gives the required nonzero coefficient.

\begin{lemma}\label{pt:subflag-restriction}
For $n\geq3$, $1\leq i\leq n-2$, and rational $m,\gamma$, one has
\begin{equation}\label{pt:three-particle-restriction}
 \left.C_n(Q;m,\gamma)\right|_{Q_i=X,\ Q_{i+1}=Y,\ Q_j=0\ (j\ne i,i+1)}
 =C_3((X,Y);m,\gamma).
\end{equation}
The substitution is understood coefficientwise in the formal
K\"ahler expansion.
\end{lemma}

\begin{proof}
Choose interlacing radii so that the open interval
$(|z_i|,|z_{i+2}|)$ contains exactly two internal vertices on row
$n-1$, one on row $n-2$, and none on the lower rows. Such a
choice is possible directly from interlacing. On row $n-1$ the
two vertices lie in the consecutive gaps on either side of
$|z_{i+1}|$. On the next row the interval between these two
vertices supplies one internal vertex; at each boundary choose
the neighboring vertex outside $(|z_i|,|z_{i+2}|)$. On the
following row all vertices can be chosen outside this interval,
and this remains possible on every lower row. Each choice is
made in a nonempty open interlacing interval, with distinct radii.

Setting $Q_j=0$ for $j\ne i,i+1$ retains only flows of degree
zero in those gaps. The sum of flows across each such cut is
zero, so no edge crossing it carries a positive index. Outside
the surviving interval all external exponents are zero; acyclicity
and internal conservation then force every remaining flow there
to vanish as well. The only active internal vertices are the
three just chosen. Together with $z_i,z_{i+1},z_{i+2}$, they
have exactly the adjacent-row edges and the single same-row edge
of the rank-three master period. All other binomial factors
contribute their constant coefficients. Relabeling the three
external variables proves \eqref{pt:three-particle-restriction}.
\end{proof}

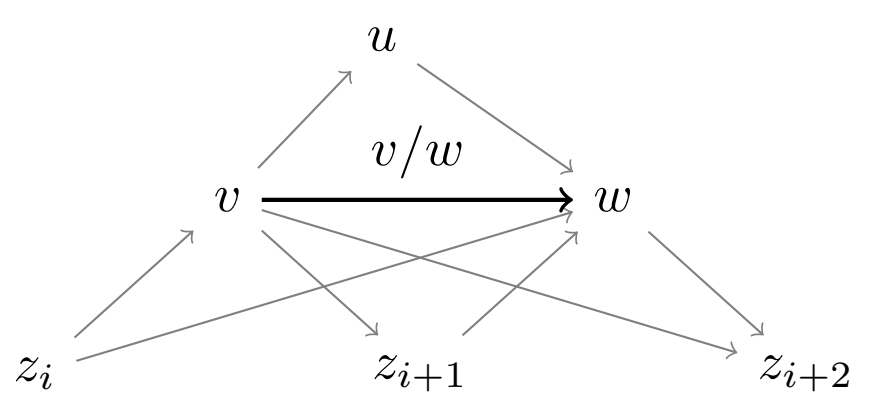
\begin{figure}[t]
\centering
\begin{tikzpicture}[scale=.85,every node/.style={font=\small}]
 \node (z1) at (0,0) {$z_i$};
 \node (z2) at (3,0) {$z_{i+1}$};
 \node (z3) at (6,0) {$z_{i+2}$};
 \node (v) at (1.5,1.35) {$v$};
 \node (w) at (4.5,1.35) {$w$};
 \node (u) at (2.7,2.6) {$u$};
 \foreach \a/\b in {z1/v,v/z2,v/z3,z1/w,z2/w,w/z3,v/u,u/w}
  {\draw[->,gray] (\a)--(\b);}
 \draw[->,thick] (v)--node[above] {$v/w$}(w);
\end{tikzpicture}
\caption{The active triangle after restricting to two consecutive
K\"ahler variables. The vertices $v,w$ lie on row $n-1$ and $u$
on row $n-2$. All edges outside this triangle have zero flow.
The bold edge is the only surviving same-row factor.}
\label{pt:subflag-figure}
\end{figure}

For the local calculation, write $u=x_{1,1}$, $v=x_{2,1}$,
$w=x_{2,2}$ and set $z_1=X$, $z_2=1$, $z_3=Y^{-1}$.
The adjacent-row product is
\begin{equation}\label{pt:edge-product}
\begin{split}
 \mathscr E_3={}&(1-v/u)(1-u/w)(1-X/v)(1-v)\\
 &\cdot(1-vY)(1-X/w)(1-1/w)(1-wY).
\end{split}
\end{equation}
Put $r=v/w$. The following weights count the flows left after
fixing the exponent of this single same-row ratio.

\begin{lemma}\label{pt:flow-weights}
For $m,h\in\Z_{\geq0}$, put
\[
 W_h^{(m)}(d,e)=[X^dY^e]\CT_{u,v,w}r^h\mathscr E_3^m.
\]
These weights are integers and vanish unless
$0\leq d,e\leq2m$ and $0\leq h\leq\min(d,e,m)$. Moreover,
\begin{equation}\label{pt:extreme-weight}
 W_h^{(m)}(h,h)=\binom mh^2.
\end{equation}
\end{lemma}

\begin{proof}
Assign indices $a,b,c,\ell,j,f,g,k$ to the eight factors of
\eqref{pt:edge-product}, in their displayed order. Each index
lies between $0$ and $m$. Conservation of the exponents of
$u,v,w$ gives
\[
 b=a,\qquad c=a+\ell+j+h,\qquad k=a+f+g+h.
\]
The external degrees are
\[
 d=c+f=a+\ell+j+f+h,\qquad
 e=j+k=a+j+f+g+h.
\]
These identities prove the support bounds. If $d=e=h$, all
indices except $c=k=h$ vanish. Their signs cancel and their
binomial weights multiply to $\binom mh^2$, proving
\eqref{pt:extreme-weight}. More generally, the sum of all eight
indices is even under the conservation laws, so every surviving
edge weight has positive sign.
\end{proof}

We can now compute the sharp coefficients in every rank.
Coefficients supported on two consecutive K\"ahler variables
are understood to have all other degrees zero.

\begin{theorem}\label{pt:sharpness}
Assume \eqref{ar:prime-assumptions}, $n\geq3$, $s\geq1$, and put
$a_0=\omega(p-1)\in\{1,\ldots,(p-1)/2\}$.
For $1\leq i\leq n-2$ and $0\leq\ell\leq s-1$,
\begin{equation}\label{pt:sharp-coefficient}
 [Q_i^{p^\ell}Q_{i+1}^{p^\ell}]\Delta_s^{(n)}
 \equiv p^{s-\ell}a_0^2\pmod{p^{s-\ell+1}}.
\end{equation}
The coefficient has valuation exactly $s-\ell$, and
\begin{equation}\label{pt:coefficient-norm}
 \|\widetilde T_s^{(n)}-T_s^{(n)}\|_p=p^{-1}
 \qquad(n\geq3,\ s\geq1).
\end{equation}
In particular,
$\widetilde T_s^{(n)}\not\equiv T_s^{(n)}\pmod{p^s}$ for every
$n\geq3$ and $s\geq2$.
\end{theorem}

\begin{proof}
Put $d_\ell=p^\ell$, an integer K\"ahler degree.
Lemma~\ref{pt:subflag-restriction}, applied to both same-row
exponents, gives
\begin{equation}\label{pt:sharp-face-coefficient}
 [Q_i^{d_\ell}Q_{i+1}^{d_\ell}]\Delta_s^{(n)}
 =[X^{d_\ell}Y^{d_\ell}]\Delta_s^{(3)}.
\end{equation}
For the right-hand side, expand the single same-row factor as
\[
 c_h=[r^h](1-r)^{-2m}\left[1-(1-r)^P\right],\qquad c_0=0.
\]
Then Lemma~\ref{pt:flow-weights} gives
\begin{equation}\label{pt:local-sharp-sum}
 [X^{d_\ell}Y^{d_\ell}]\Delta_s^{(3)}
 =\sum_{h=1}^{d_\ell}c_hW_h^{(m)}(d_\ell,d_\ell).
\end{equation}
Since $m=a_0(1+p+\cdots+p^{s-1})$ and $a_0\geq1$, we have
$d_\ell\leq m<P$.

If $\ell=0$, then $d_\ell=1$, $c_1=P$, and
$W_1^{(m)}(1,1)=m^2$. The coefficient is therefore $Pm^2$,
which proves \eqref{pt:sharp-coefficient} because
$m\equiv a_0\pmod p$.
For $\ell\geq1$, every $1\leq k<d_\ell=p^\ell$ has
$v_p(k)\leq\ell-1$. Equation~\eqref{pt:binomial-valuation}
therefore gives
$c_h\equiv0\pmod{p^{s-\ell+1}}$ for $h<d_\ell$.
In $c_{d_\ell}$, all terms with binomial index $k<d_\ell$
vanish at the same precision. The remaining term has positive sign because
$d_\ell$ is odd, so
\[
 c_{d_\ell}\equiv\binom P{d_\ell}
 =p^{s-\ell}\binom{P-1}{d_\ell-1}
 \equiv p^{s-\ell}\pmod{p^{s-\ell+1}}.
\]
Here $\binom{P-1}{d_\ell-1}\equiv(-1)^{d_\ell-1}=1\pmod p$.
Using \eqref{pt:extreme-weight} in
\eqref{pt:local-sharp-sum}, we obtain
\[
 [X^{d_\ell}Y^{d_\ell}]\Delta_s^{(3)}
 \equiv p^{s-\ell}\binom m{d_\ell}^2\pmod{p^{s-\ell+1}}.
\]

The base-$p$ digits of $m$ are all $a_0$. Equivalently,
\[
 (1+X)^m\equiv\prod_{j=0}^{s-1}(1+X^{p^j})^{a_0}\pmod p.
\]
The factors with $j<\ell$ have total degree
$a_0(d_\ell-1)/(p-1)\leq(d_\ell-1)/2<d_\ell$.
Those with $j>\ell$ have no nonconstant term of degree at most
$d_\ell$. Thus the
coefficient of $X^{d_\ell}$ is $a_0$, and
$\binom m{d_\ell}\equiv a_0\pmod p$. This proves
\eqref{pt:sharp-coefficient} for every $i$ and $n$.

Theorem~\ref{pt:coefficientwise-bound} gives the norm bound
$\|\Delta_s^{(n)}\|_p\leq p^{-1}$. Taking $\ell=s-1$ in
\eqref{pt:sharp-coefficient} gives a coefficient of valuation
exactly one, proving equality in \eqref{pt:coefficient-norm}.
For $s\geq2$, this coefficient also excludes divisibility by
$p^s$.
\end{proof}

Thus every complete flag of rank at least three exhibits the
same distinction between the two modes of approximation. The
coefficients converge to the same characteristic-zero partition
function at every fixed degree, but the difference retains norm
$p^{-1}$. The sharp coefficients lie on every consecutive
three-particle subflag. The higher-precision Dwork congruences
proved in this paper concern the original family $T_s^{(n)}$;
the comparison does not assert these congruences for
$\widetilde T_s^{(n)}$.

\begin{corollary}\label{pt:distance-to-limit}
For every admissible $\omega,p$, every $n\geq2$ and $s\geq1$,
\[
 \|T_s^{(n)}-V_n(Q;\omega)\|_p
 =\|\widetilde T_s^{(n)}-V_n(Q;\omega)\|_p=1.
\]
Thus neither family converges to the vertex in coefficient norm.
\end{corollary}

\begin{proof}
Restricting the normalized quadratic equation to a simple-root
axis gives the two-particle equation. Characteristic-zero
uniqueness therefore identifies the restrictions of $V_n$ and
$F_m^{(n)}$ with $V_2$ and $F_m^{(2)}$, respectively.
For $n\geq3$, Theorem~\ref{pt:coefficientwise-bound} gives
the same axis restriction for $T_s^{(n)}$ as for
$\widetilde T_s^{(n)}$; for $n=2$ they coincide directly.
With $P=p^s$ and $m=A_s<P$, both polynomial restrictions have
degree $m$, whereas
\[
 [Q_i^P]V_n=\binom{-\omega}{P}^{\!2}\equiv a_0^2\pmod p,
 \qquad a_0=\omega(p-1)\ne0\pmod p.
\]
Indeed $-\omega=a_0(1+p+p^2+\cdots)$ in $\Z_p$.
For every $t>s$, the base-$p$ digits of $A_t$ equal $a_0$,
so Lucas congruences give $\binom{A_t}{p^s}\equiv a_0\pmod p$.
Taking the coefficientwise limit proves the displayed congruence.
Each difference has integral coefficients, and its coefficient
at $Q_i^P$ is a unit. Its norm is consequently exactly one.
\end{proof}

We next express the comparison in the monomial symmetric
basis. This avoids specializing a
full Jack basis at a negative rational parameter.

\begin{corollary}\label{pt:symmetric-comparison}
Assume \eqref{ar:prime-assumptions} and $n\geq3$, put $m=A_s$ and
$\lambda=m\delta$, and use $\mathcal T_s$ from \eqref{pt:symmetric-master}.
Write $m_\mu$ for the sum of the distinct
monomials in the permutation orbit of $z^\mu$, and expand
\begin{equation}\label{pt:symmetric-discrepancy}
 P_\lambda^{(-1/m)}-\mathcal T_s
 =\sum_{\mu\leq\lambda}c_s(\mu)m_\mu.
\end{equation}
Here the sum is over partitions of $|\lambda|$ of length at
most $n$, ordered by dominance, and $c_s(\mu)\in\Z$.
For such a partition, padded by zeros to length $n$, put
\[
 \sigma_j(\mu)
 =\sum_{a=n-j+1}^n\mu_a-\frac{mj(j-1)}2,
 \qquad
 b_m(\mu)=\max_{1\leq j\leq n-2}
                 \min\{\sigma_j(\mu),\sigma_{j+1}(\mu)\}.
\]
If $b_m(\mu)=0$, then $c_s(\mu)=0$. Otherwise,
\begin{equation}\label{pt:symmetric-precision}
 v_p(c_s(\mu))\geq
 \max\{1,\ s-\lfloor\log_pb_m(\mu)\rfloor\}.
\end{equation}
More explicitly, for $1\leq k\leq n-2$ and $0\leq\ell\leq s-1$,
\begin{equation}\label{pt:three-row-squeeze}
 c_s\bigl(\lambda-p^\ell\varepsilon_k+p^\ell\varepsilon_{k+2}\bigr)
 \equiv p^{s-\ell}a_0^2\pmod{p^{s-\ell+1}},
 \qquad a_0=\omega(p-1).
\end{equation}
Consequently, for every $n\geq3$ and $s\geq1$,
\[
 P_{A_s\delta}^{(-1/A_s)}\equiv\mathcal T_s\pmod p,
 \qquad
 P_{A_s\delta}^{(-1/A_s)}\not\equiv\mathcal T_s\pmod{p^2}.
\]
\end{corollary}

\begin{proof}
The dependence of $\Phi_s(x,z)$ on the external variables is
through $E(x,z)^m$. Proposition~\ref{pt:positive-support}
gives the symmetry, integrality and common permutahedral Newton
polytope
of the two sides in \eqref{pt:symmetric-discrepancy}.

To identify their coefficients, choose the increasing exponent
vector $\beta=(\mu_n,\ldots,\mu_1)$ in the orbit of $\mu$.
Dividing $z^\beta$ by $G_n^m$ gives the chamber monomial
$Q^{\boldsymbol\sigma(\mu)}$. Symmetry and
\eqref{pt:jack-specialization} consequently give
\[
 c_s(\mu)
 =[Q^{\boldsymbol\sigma(\mu)}]\Delta_s^{(n)}.
\]
Dominance makes every $\sigma_j(\mu)$ nonnegative.
Theorem~\ref{pt:coefficientwise-bound} proves the vanishing
assertion and \eqref{pt:symmetric-precision}.

For \eqref{pt:three-row-squeeze}, put $d_\ell=p^\ell$.
The inequality $d_\ell\leq m$ shows that
$\mu=\lambda-d_\ell\varepsilon_k+d_\ell\varepsilon_{k+2}$
is still a partition. Reversing its entries gives
\[
 \beta=(0,m,\ldots,(n-1)m)
       +d_\ell\varepsilon_i-d_\ell\varepsilon_{i+2},
 \qquad i=n-k-1.
\]
Thus $\sigma_i(\mu)=\sigma_{i+1}(\mu)=d_\ell$, and all its
other chamber degrees vanish.
Theorem~\ref{pt:sharpness} gives the stated residue.
Taking $\ell=s-1$ yields valuation exactly one, proving the
last assertion. Since $m_\mu$ sums distinct monomials, this
argument introduces no orbit-size or factorial denominator.
\end{proof}

The degree estimate also determines the discrepancy in weighted
coefficient norms, which measure convergence on closed
polydiscs of different radii.

\begin{corollary}\label{pt:weighted-gauss}
Assume \eqref{ar:prime-assumptions}, $n\geq3$, and $s\geq1$.
For a real number $0<\varrho\leq1$, define
\[
 \left\|\sum_{\boldsymbol e}c_{\boldsymbol e}Q^{\boldsymbol e}
 \right\|_{p,\varrho}
 =\sup_{\boldsymbol e}|c_{\boldsymbol e}|_p
      \varrho^{|\boldsymbol e|},
 \qquad |\boldsymbol e|=e_1+\cdots+e_{n-1}.
\]
Then
\begin{equation}\label{pt:weighted-gauss-distance}
 \|\Delta_s^{(n)}\|_{p,\varrho}
 =\max_{0\leq\ell\leq s-1}
      p^{-(s-\ell)}\varrho^{2p^\ell}.
\end{equation}
At $\varrho=1$ this is $p^{-1}$. For $0<\varrho<1$, put
\[
 C_{\varrho}=\max_{\ell\geq0}p^\ell\varrho^{2p^\ell}<\infty.
\]
For every $s$ one has
$\|\Delta_s^{(n)}\|_{p,\varrho}\leq C_{\varrho}p^{-s}$,
and equality holds for all sufficiently large $s$.
In particular, the discrepancy converges to zero in every such
weighted norm, with rate $O_{\varrho}(p^{-s})$.
\end{corollary}

\begin{proof}
Let $c_{\boldsymbol e}=[Q^{\boldsymbol e}]\Delta_s^{(n)}$.
If $B(\boldsymbol e)=0$, Theorem~\ref{pt:coefficientwise-bound}
gives $c_{\boldsymbol e}=0$. Otherwise set
\[
 \ell=\min\{s-1,\lfloor\log_p B(\boldsymbol e)\rfloor\}.
\]
The same theorem gives $v_p(c_{\boldsymbol e})\geq s-\ell$.
Since some adjacent degrees are both at least
$B(\boldsymbol e)$,
$|\boldsymbol e|\geq2B(\boldsymbol e)\geq2p^\ell$.
Consequently
\[
 |c_{\boldsymbol e}|_p\varrho^{|\boldsymbol e|}
 \leq p^{-(s-\ell)}\varrho^{2p^\ell},
\]
which proves the upper bound in
\eqref{pt:weighted-gauss-distance}.
For every $0\leq\ell\leq s-1$, Theorem~\ref{pt:sharpness}
provides a coefficient with precisely two adjacent degrees
$p^\ell,p^\ell$ and valuation $s-\ell$. Its weighted absolute
value attains the corresponding term of the maximum, proving
equality.

If $\varrho<1$, the sequence
$p^\ell\varrho^{2p^\ell}$ tends to zero. Its positive maximum
is therefore finite and attained at a finite index. Once this
index lies below $s$, formula~\eqref{pt:weighted-gauss-distance}
equals $C_{\varrho}p^{-s}$, proving the remaining assertions.
\end{proof}

\section{Dwork Congruences and Differential Equations}
\label{ar:section-ghosts}
We prove the Dwork congruences for the master-polynomial family,
then use their Frobenius quotients to reconstruct the vertex by
an infinite product. A differential consequence of these
congruences transfers the full commuting tCM equations from
the vertex to its polynomial approximations.

\subsection{Alternating Ghosts}
The master polynomials have two structures adapted to the
Frobenius action: their exponents grow with powers of $p$, and
they are alternating in each internal row. We shall use these
properties to prove the Dwork congruences for $T_s^{(n)}$.
Subtracting products obtained from lower levels by Frobenius
substitution defines the ghost polynomial $K_s$. Its coefficients
are divisible by $p^{s-1}$. Degree bounds and alternation then
determine all coefficients relevant to extraction from a single
choice of exponents in each row. The polynomial identity for
$K_s$ consequently becomes a recurrence for $T_s^{(n)}$ alone.

We follow the ghost method originating in Mellit's approach to
Dwork congruences, developed for Laurent-polynomial constant terms
in \cite{MV16} (\S2) and for hypergeometric master polynomials in
\cite{VZ24} (\S2). Its alternating form appears in \cite{SV23}
(\S\S5.2--5.3). For the general relation between Frobenius and
$p$-adic differential equations we refer to \cite{Ked10}.

Let $\sigma$ denote the substitution $(x,z)\mapsto(x^p,z^p)$ on
integral polynomials, acting trivially on coefficients. Its action
on K\"ahler series is $Q\mapsto Q^p$. The elementary ghost terms
associated with an integral polynomial $f$ are
\[
 R_j(f)=f^{p^j}-\sigma\bigl(f^{p^{j-1}}\bigr)
       \in p^j\Z[x,z],\qquad j\geq1.
\]
For $j=1$ this is $f^p\equiv\sigma(f)\pmod p$.
Raising congruent polynomials to successive $p$th powers gives the
higher divisibility. The ghost terms therefore measure the error
in replacing ordinary powers by Frobenius substitution at each
$p$-adic precision.

For the master polynomials, the level-$s$ exponents contain the
sum $1+p+\cdots+p^{s-1}$. We separate its successive powers of
$p$ by the following factorization and recursion.

Put
\begin{equation}\label{ar:bar-Phi}
 f=M^{\beta_1}\Delta^{\gamma_1-1}E^{A_1},\qquad
 \mathcal B_0=1,\qquad
 \mathcal B_s=f^{1+p+\cdots+p^{s-1}}\quad(s\geq1).
\end{equation}
Recall that $\mathfrak S=\prod_{i=1}^{n-1}S_i$ permutes the
variables independently in each internal row and fixes the
external variables. For $\tau=(\tau_1,\ldots,\tau_{n-1})$, put
$\operatorname{sgn}(\tau)=\prod_i\operatorname{sgn}(\tau_i)$.
A polynomial is invariant if every row permutation leaves it
unchanged, and alternating if each permutation multiplies it by
its sign.
We have $\Phi_s=\Delta\mathcal B_s$. The polynomials $f$ and
$\mathcal B_s$ are invariant, while $\Delta$ is alternating.

\begin{lemma}\label{ar:ghost-divisibility}
Define $K_s$, $s\geq1$, recursively by
\begin{equation}\label{ar:K-recursion}
 K_s=\Phi_s-\sum_{j=1}^{s-1}K_j\sigma^j(\mathcal B_{s-j}).
\end{equation}
Then $K_s$ is alternating under $\mathfrak S$ and
\[
 K_s\in p^{s-1}\Z[x,z].
\]
For every variable in row $i$,
\begin{equation}\label{ar:degree-bound}
 \deg_{x_{i,a}}K_s\leq D_{i,s},\qquad
 D_{i,s}=(i+\omega)(p^s-1)+i-1<(i+1)p^s-1.
\end{equation}
The same bound holds for $\Phi_s$.
\end{lemma}

\begin{proof}
Write $K_s=\Delta\ell_s$, which is possible inductively from the
recursion. Dividing the polynomial identity by its common factor
$\Delta$ gives
\[
 \mathcal B_s=\sum_{j=1}^s\ell_j\sigma^j(\mathcal B_{s-j}).
\]
Induction shows that $\ell_s$ is symmetric, so $K_s$ is
alternating. Using $\mathcal B_{s+1}=\mathcal B_sf^{p^s}$, we
subtract the recursions at consecutive levels and obtain
\begin{equation}\label{ar:ghost-power-identity}
 \ell_{s+1}=
 \sum_{j=1}^s\ell_j\sigma^j(\mathcal B_{s-j})
 \left(f^{p^s}-\sigma^j(f^{p^{s-j}})\right).
\end{equation}
For integral polynomials $u\equiv v\pmod p$, induction using the
binomial theorem gives
$u^{p^k}\equiv v^{p^k}\pmod{p^{k+1}}$. Applying this to
$f^{p^j}\equiv\sigma^j(f)\pmod p$ yields
\[
 f^{p^s}-\sigma^j(f^{p^{s-j}})\in p^{s-j+1}\Z[x,z].
\]
Since $\ell_1=f$, induction in \eqref{ar:ghost-power-identity}
proves $\ell_s\in p^{s-1}\Z[x,z]$ and the required divisibility
of $K_s$. We have used the factorization by $\Delta$ only in
polynomial identities, without dividing a congruence by it.

To estimate the degrees, observe that a variable in row $i$
occurs in $i-1$ same-row factors and $2i$ adjacent-row or framing
factors. Hence
\[
 \deg_{x_{i,a}}\Phi_s
 \leq\beta_s+(i-1)\gamma_s+2iA_s=D_{i,s},
 \qquad
 \deg_{x_{i,a}}\mathcal B_s\leq(i+\omega)(p^s-1).
\]
Each product in \eqref{ar:K-recursion} has degree bounded by
\[
 D_{i,j}+p^j(i+\omega)(p^{s-j}-1)=D_{i,s},
\]
and induction gives the bound for $K_s$. The strict inequality
follows from
$(i+1)p^s-1-D_{i,s}=(1-\omega)p^s+\omega>0$.
\end{proof}

The degree estimate determines which coefficients can contribute
to the Frobenius recursion. An exponent congruent to $-1$ modulo
$p^s$ in row $i$ must be one of $p^s-1,\ldots,ip^s-1$.
Alternation excludes repetitions, so the remaining choices are
permutations of these exponents.

\begin{lemma}\label{ar:orbit}
For $F_s=\Phi_s$ or $F_s=K_s$, a coefficient
$[x^{p^s\boldsymbol u-\mathbf1}]F_s$ with
$\boldsymbol u\in\Z^N$ can be nonzero only when
$\boldsymbol u=\tau(\mathbf d)$ for $\tau\in\mathfrak S$. In that case
\[
 [x^{p^s\tau(\mathbf d)-\mathbf1}]F_s
 =\operatorname{sgn}(\tau)
   [x^{p^s\mathbf d-\mathbf1}]F_s.
\]
\end{lemma}

\begin{proof}
By nonnegativity and \eqref{ar:degree-bound},
$1\leq u_{i,a}\leq i$. If two entries of row $i$ coincide,
their transposition fixes the monomial and changes the sign of
its coefficient. This coefficient vanishes over $\Z$.
The entries are therefore distinct and form a permutation of
$1,\ldots,i$. Alternation gives the corresponding sign.
\end{proof}

Consider $K_j\sigma^j(\mathcal B_{s-j})$. Since the exponents
of the second factor are divisible by $p^j$, only the coefficients
of $K_j$ in Lemma~\ref{ar:orbit} can contribute to our extraction.
Summing these coefficients with their permutation signs
reproduces the Vandermonde factor in the lower-level
coefficient $U_{s-j}$. The resulting recurrence involves only
the K\"ahler variables.

\begin{proposition}\label{ar:scalar-ghosts}
There are polynomials $H_s\in\Z[Q_1,\ldots,Q_{n-1}]$, $s\geq1$, such
that
\begin{equation}\label{ar:T-ghost-recursion}
 T_s^{(n)}=\sum_{j=1}^s H_j\sigma^j(T_{s-j}^{(n)}),
 \qquad H_s\in p^{s-1}\Z[Q_1,\ldots,Q_{n-1}].
\end{equation}
\end{proposition}

\begin{proof}
Let $\Gamma_j(z)=[x^{p^j\mathbf d-\mathbf1}]K_j(x,z)$.
Expanding the Vandermonde determinant gives
\begin{equation}\label{ar:determinant-expansion}
 \Delta(x)=\sum_{\tau\in\mathfrak S}
       \operatorname{sgn}(\tau)x^{\tau(\mathbf d)-\mathbf1}.
\end{equation}
It follows that, for every $l\geq0$,
\begin{equation}\label{ar:determinant-extraction}
 U_l(z)=\sum_{\tau\in\mathfrak S}\operatorname{sgn}(\tau)
 [x^{p^l\mathbf d-\tau(\mathbf d)}]\mathcal B_l(x,z).
\end{equation}

We extract $x^{p^s\mathbf d-\mathbf1}$ from
$K_j\sigma^j(\mathcal B_{s-j})$. The contributing exponents
satisfy $\alpha+p^j\beta=p^s\mathbf d-\mathbf1$, so
$\alpha=p^j\boldsymbol u-\mathbf1$ for an integral vector
$\boldsymbol u$. Lemma~\ref{ar:orbit} gives
$\boldsymbol u=\tau(\mathbf d)$ and coefficient
$\operatorname{sgn}(\tau)\Gamma_j$. The exponent in the second
factor is $\beta=p^{s-j}\mathbf d-\tau(\mathbf d)$.
Summing over $\tau$ and using \eqref{ar:determinant-extraction},
we find
\[
 [x^{p^s\mathbf d-\mathbf1}]
 K_j\sigma^j(\mathcal B_{s-j})
 =\Gamma_j(z)U_{s-j}(z^{p^j}).
\]
The identity also holds for $j=s$, since $\mathcal B_0=1$ and
$U_0=1$. Therefore
\[
 U_s=\sum_{j=1}^s\Gamma_j\sigma^j(U_{s-j}).
\]

In $\Z[z_1^{\pm1},\ldots,z_n^{\pm1}]$, set
\[
 H_j=(-1)^{A_jc_n}G_n^{-A_j}\Gamma_j.
\]
The identities $A_s=A_j+p^jA_{s-j}$ and
$G_n(z^{p^j})=G_n(z)^{p^j}$ match the monomial normalizations.
Since $p$ is odd, the signs satisfy
$(-1)^{A_sc_n}=(-1)^{A_jc_n}(-1)^{A_{s-j}c_n}$ as well.
Normalizing the recurrence for $U_s$ thus gives
\eqref{ar:T-ghost-recursion}. This recurrence and polynomiality
of $T_s^{(n)}$ imply inductively that $H_s$ is a polynomial in
$Q$. Its Laurent coefficients are divisible by $p^{s-1}$ by
Lemma~\ref{ar:ghost-divisibility}. Distinct $Q$-monomials give
distinct $z$-Laurent monomials, proving the same divisibility in
$\Z[Q]$.
\end{proof}

The recurrence no longer contains internal variables. Its
divisibility property gives congruences between successive
Frobenius quotients by the following formal argument, applicable
to any series with constant coefficient one.

\begin{lemma}\label{ar:formal-ghost}
Suppose $F_0=1$ and $F_s\in1+I R$, where
$R=\Z_p[[Q_1,\ldots,Q_{n-1}]]$ and $I=(Q_1,\ldots,Q_{n-1})$.
If
\[
 F_s=\sum_{j=1}^s H_j\sigma^j(F_{s-j}),\qquad H_j\in p^{j-1}R,
\]
then the formal units $\mathcal R_s=F_s/\sigma(F_{s-1})$ satisfy
\begin{equation}\label{ar:ratio-congruence}
 \mathcal R_{s+1}-\mathcal R_s\in p^sR\qquad(s\geq1).
\end{equation}
\end{lemma}

\begin{proof}
For $s=1$ we use $F_2=H_1\sigma(F_1)+H_2$ and $F_1=H_1$.
Suppose the assertion holds for all indices smaller than $s$.
Dividing the ghost recurrences by their Frobenius-shifted
denominators gives
\[
 \mathcal R_{s+1}=
 \sum_{j=1}^{s+1}H_j
   \frac{\sigma^j(F_{s+1-j})}{\sigma(F_s)},\qquad
 \mathcal R_s=
 \sum_{j=1}^{s}H_j
   \frac{\sigma^j(F_{s-j})}{\sigma(F_{s-1})}.
\]
The terms with $j=1$ coincide, and the last term of the first sum
belongs to $p^sR$. For $2\leq j\leq s$ we use the telescoping
identities
\[
 \frac{F_s}{\sigma^{j-1}(F_{s-j+1})}
 =\prod_{i=0}^{j-2}\sigma^i(\mathcal R_{s-i}),\qquad
 \frac{F_{s-1}}{\sigma^{j-1}(F_{s-j})}
 =\prod_{i=0}^{j-2}\sigma^i(\mathcal R_{s-i-1}).
\]
The induction hypothesis identifies these products modulo
$p^{s-j+1}$. Their factors are units, so we may invert and apply
$\sigma$. Multiplication by $H_j\in p^{j-1}R$ then places the
difference of the corresponding summands in $p^sR$, proving
\eqref{ar:ratio-congruence}.
\end{proof}

For the master-polynomial recurrence, the condition
$T_s^{(n)}(0)=1$ makes all these denominators invertible. We can
therefore clear them and express the ratio congruence as a
polynomial identity modulo $p^s$.

\begin{theorem}\label{ar:dwork}
Under \eqref{ar:prime-assumptions}, the polynomials
$T_s^{(n)}$ satisfy, for every $s\geq1$,
\begin{equation}\label{ar:dwork-cross}
 T_{s+1}^{(n)}(Q)T_{s-1}^{(n)}(Q^p)
 -T_s^{(n)}(Q)T_s^{(n)}(Q^p)
 \in p^s\Z[Q_1,\ldots,Q_{n-1}].
\end{equation}
Equivalently, their successive formal Frobenius ratios are congruent
modulo $p^s$ as in \eqref{ar:ratio-congruence}.
\end{theorem}

\begin{proof}
The recurrence of Proposition~\ref{ar:scalar-ghosts} satisfies
Lemma~\ref{ar:formal-ghost}; its denominators are units by
Proposition~\ref{ar:normalization}. Cross multiplication gives
the congruence in $\Z_p[[Q]]$. Since the left-hand side is an
integral polynomial, its coefficients are divisible by $p^s$
in $\Z[Q]$ as claimed.
\end{proof}

\begin{remark}\label{ar:small-primes}
The proof does not divide by $|\mathfrak S|=\prod_{i=1}^{n-1}i!$.
Consequently the congruences remain valid for $p\leq n$.
\end{remark}

\subsection{Infinite Products for the Vertex}
The Frobenius-quotient construction follows Dwork's approach to
period congruences \cite{D69}.
\label{ar:section-product}

The Dwork congruences reconstruct the partition function as an
infinite product. They give a limit of the ratios
$T_s^{(n)}(Q)/T_{s-1}^{(n)}(Q^p)$, which we identify with
$\mathcal V_n(Q)/\mathcal V_n(Q^p)$. Iterating this identity
reconstructs the period from its Frobenius quotient. Modulo
$p^m$, all factors involve the same pair $T_m^{(n)}$ and
$T_{m-1}^{(n)}$. Their nonconstant terms have increasing order
under $Q\mapsto Q^p$, so the product exists as a formal series.

\begin{theorem}\label{ar:product}
Under \eqref{ar:prime-assumptions}, for every integer $m\geq1$,
\begin{equation}\label{ar:infinite-product}
 \mathcal V_n(Q;\omega)\equiv
 \prod_{i=0}^{\infty}
 \frac{T_m^{(n)}(Q^{p^i})}
      {T_{m-1}^{(n)}(Q^{p^{i+1}})}
 \pmod{p^m}.
\end{equation}
The product is taken in the $I$-adic topology for
$I=(Q_1,\ldots,Q_{n-1})$. Retaining its first $K$ factors gives
an error in
\begin{equation}\label{ar:product-cutoff}
 p^m\Z_p[[Q]]+I^{p^K}.
\end{equation}
In particular, coefficients of total degree at most $D$ require
only $K$ factors with $p^K>D$.
\end{theorem}

\begin{proof}
Write $T_s=T_s^{(n)}$, $\mathcal V=\mathcal V_n$, and
$R=\Z_p[[Q]]$, and put
$\mathcal R_s=T_s/\sigma(T_{s-1})$. By
Theorem~\ref{ar:dwork}, these ratios form a uniformly $p$-adically
Cauchy sequence. The completeness of $R$ gives a limit
$\mathcal R$, and telescoping the consecutive congruences yields
\begin{equation}\label{ar:limit-ratio-precision}
 \mathcal R-\mathcal R_m\in p^mR.
\end{equation}
For a series with constant coefficient one, every coefficient of
its inverse is a polynomial in finitely many coefficients of the
series. Formal inversion is therefore continuous in the
coefficientwise $p$-adic topology. Applying
Theorem~\ref{ar:limit}, we obtain
$\mathcal R=\mathcal V/\sigma(\mathcal V)$.

The finite products satisfy
\[
 \prod_{i=0}^{K-1}\sigma^i(\mathcal R)
 =\frac{\mathcal V}{\sigma^K(\mathcal V)}.
\]
Since $\mathcal V\in1+I$, the denominator belongs to
$1+I^{p^K}$. Equation \eqref{ar:limit-ratio-precision} allows
us to replace each factor by $\sigma^i(\mathcal R_m)$ modulo
$p^mR$, giving \eqref{ar:product-cutoff}. Moreover,
$\sigma^i(\mathcal R_m)\in1+I^{p^i}$, so the products converge
$I$-adically. Taking this limit proves
\eqref{ar:infinite-product}.
\end{proof}

The same calculation applies before taking the period limit.
Stopping the telescoping product at a finite level leaves a
lower-level polynomial as its last factor.

\begin{remark}\label{ar:finite-product}
For $s\geq m$ we have
\begin{equation}\label{ar:finite-level-product}
 T_s^{(n)}\equiv
 \left(\prod_{i=0}^{s-m}
   \frac{T_m^{(n)}(Q^{p^i})}
        {T_{m-1}^{(n)}(Q^{p^{i+1}})}\right)
 T_{m-1}^{(n)}(Q^{p^{s-m+1}})\pmod{p^m}.
\end{equation}
Indeed, telescope
$T_s=\prod_{i=0}^{s-1}\sigma^i(\mathcal R_{s-i})$ and replace
the ratios with index at least $m$ modulo $p^m$. Recombining the
remaining factors gives the stated formula. Dwork congruences
and constant coefficient one therefore imply a coefficientwise
$p$-adic limit, identified with the master period by
Theorem~\ref{ar:limit}. The infinite product converges in the
$I$-adic topology; the exact distance in Corollary~\ref{pt:distance-to-limit} shows
that the polynomials themselves do not converge in coefficient norm.
\end{remark}

Modulo $p$, the identity $T_0^{(n)}=1$ removes the denominator
factors. We obtain product formulas involving only
Frobenius substitutions in $T_1^{(n)}$.

\begin{corollary}\label{ar:mod-p-product}
One has
\[
 \mathcal V_n(Q;\omega)\equiv
 \prod_{i=0}^{\infty}T_1^{(n)}(Q^{p^i})\pmod p,
 \qquad
 T_s^{(n)}\equiv\prod_{i=0}^{s-1}T_1^{(n)}(Q^{p^i})\pmod p.
\]
\end{corollary}

\begin{proof}
Set $m=1$ in Theorem~\ref{ar:product} and
\eqref{ar:finite-level-product}, respectively.
\end{proof}

Over $\F_p$, Frobenius substitution agrees with taking the
$p$th power of a series. Thus the product relation also gives an
algebraic equation for the reduced period.

\begin{corollary}\label{ar:algebraicity}
The reduced period satisfies
\[
 T_1^{(n)}(Q)\mathcal V_n(Q;\omega)^{p-1}=1
 \quad\text{in }\F_p[[Q_1,\ldots,Q_{n-1}]].
\]
It is the unique solution of this equation with constant coefficient one.
\end{corollary}
\begin{proof}
The product formula gives
$\mathcal V_n(Q)=T_1^{(n)}(Q)\mathcal V_n(Q^p)$ modulo $p$.
Using $\mathcal V_n(Q^p)=\mathcal V_n(Q)^p$ in characteristic
$p$ and dividing by the unit $\mathcal V_n$, we obtain the
equation. At $(Q,Y)=(0,1)$ the derivative of
$T_1^{(n)}Y^{p-1}-1$ with respect to $Y$ is $p-1\ne0$;
the coefficient recursion therefore gives uniqueness.
\end{proof}

Theorem~\ref{con:vertex} identifies $\mathcal V_n=V_n$ in all these formulas.

\subsection{The Full Commuting Congruences}
The Dwork congruences and the integral determinant coefficients
give a short proof of the finite-level equations. We first
isolate the differential consequence of the congruences; it
applies to every integral differential annihilator, without a
commutativity assumption. The integral operators of Section~\ref{cg:section} then give
the joint tCM congruences without a restriction on differential order.

\begin{lemma}\label{tc:dwork-transfer}
Let $R=\Z_p[[Q_1,\ldots,Q_m]]$, let
$\sigma f(Q)=f(Q^p)$, and put $\theta_i=Q_i\partial_{Q_i}$.
Suppose $T_s\in1+(Q_1,\ldots,Q_m)R$, $s\geq0$, satisfy
$T_0=1$, converge coefficientwise $p$-adically to $V$, and obey
\[
 T_{s+1}\sigma(T_{s-1})\equiv T_s\sigma(T_s)\pmod{p^s},
 \qquad s\geq1.
\]
For $E_s=T_s/V$ and every nonzero multi-index $\alpha$,
\begin{equation}\label{tc:dwork-flatness}
 \theta^\alpha E_s\in p^sR.
\end{equation}
If $L=\sum_\alpha a_\alpha(Q)\theta^\alpha$ is any finite-order
differential operator with $a_\alpha\in R$, then
\begin{equation}\label{tc:dwork-transfer-equation}
 LT_s\equiv E_sLV\pmod{p^s}.
\end{equation}
In particular, $LV=0$ implies $LT_s\equiv0\pmod{p^s}$.
\end{lemma}

\begin{proof}
Set $\mathcal R_s=T_s/\sigma(T_{s-1})$ for $s\geq1$.
The Dwork congruences give
$\mathcal R_{s+1}-\mathcal R_s\in p^sR$, so completeness gives
a limit $\mathcal R$ with
$\mathcal R_s-\mathcal R\in p^sR$.
Coefficientwise convergence of $T_s$, and continuity of formal
inversion for series with constant coefficient one, identify
$\mathcal R=V/\sigma(V)$. Hence
\[
 E_s-\sigma(E_{s-1})
 =\frac{\sigma(T_{s-1})}{V}
   \left(\mathcal R_s-\frac{V}{\sigma(V)}\right)
 \in p^sR.
\]
Every Euler derivative preserves $p^sR$. Since
$\theta^\alpha\sigma=p^{|\alpha|}\sigma\theta^\alpha$,
induction on $s$, starting with the integral series $E_0=1/V$,
proves \eqref{tc:dwork-flatness}. In Leibniz's rule for
$L(E_sV)$, every term except $E_sLV$ differentiates $E_s$
at least once. All these terms are divisible by $p^s$,
proving \eqref{tc:dwork-transfer-equation}. No integer or
differential-order factor is inverted.
\end{proof}

\begin{theorem}\label{tc:full-congruence}\label{tc:quadratic-congruence}
Assume \eqref{ar:prime-assumptions}. For every $n\geq2$,
$s\geq1$, and $1\leq k\leq n$,
\begin{equation}\label{tc:full-congruence-system}
 \mathsf I_{n,k}(J_nT_s^{(n)})\equiv0\pmod{p^s},\qquad
 \mathsf B_k^{(0)}T_s^{(n)}\equiv0\pmod{p^s}
 \quad\text{in }\Z_p[[Q]].
\end{equation}
Equivalently, the generating determinant satisfies
\begin{equation}\label{tc:generating-congruence}
 \mathscr S_n(u)T_s^{(n)}\equiv u^nT_s^{(n)}\pmod{p^s}.
\end{equation}
The center equation is exact, and no restriction $p>n$ is needed.
\end{theorem}

\begin{proof}
Theorems~\ref{ar:dwork}, \ref{ar:limit}, and~\ref{con:vertex}
place $T_s=T_s^{(n)}$ and $V=V_n$ under the hypotheses of
Lemma~\ref{tc:dwork-transfer}. The coefficients of
$\mathsf B_k^{(0)}$ are integral by
Lemma~\ref{tc:integral-operators}, and
$\mathsf B_k^{(0)}V_n=0$ by
Proposition~\ref{dg:full-family}. On ratio series,
$D_i=\theta_i-\theta_{i-1}$, so the transfer lemma gives
$\mathsf B_k^{(0)}T_s^{(n)}\equiv0\pmod{p^s}$.
Conjugation by the integral unit $J_n$ proves the
Schr\"odinger congruence, and \eqref{tc:determinant} gives
\eqref{tc:generating-congruence}. The center equation is exact
because $J_nT_s^{(n)}$ depends only on ratios.

\end{proof}

Every polynomial with $\Z_p$ coefficients and zero constant term
in the commuting operators $\mathsf I_{n,k}$ satisfies the same
congruence on $J_nT_s^{(n)}$. This includes the integral Newton
power-sum Hamiltonians: Newton's recursion uses integer
coefficients and multiplication by $k$, with no division by $k$.
For $k=2$, the relation
$\mathsf I_{n,2}=\frac12(\sum_iD_i)^2-H_{n,\omega}$ recovers
the quadratic congruence
$H_{n,\omega}(J_nT_s^{(n)})\equiv0\pmod{p^s}$.

\begin{remark}\label{pt:frobenius-ambiguity}
Fix $n\geq3$ and $1\leq i\leq n-2$. Put $Z=Q_iQ_{i+1}$ and
\[
 U_s(Q)=1+pZ^{p^{s-1}},\qquad
 S_s^{(n)}(Q)=T_s^{(n)}(Q)U_s(Q)\quad(s\geq1),\qquad S_0^{(n)}=1.
\]
These polynomials satisfy the full joint equations modulo $p^s$
and the Dwork congruences, converge coefficientwise to $V_n$,
and obey $\|S_s^{(n)}-T_s^{(n)}\|_p=p^{-1}$.
Indeed, every positive-order Euler derivative of
$pZ^{p^{s-1}}$ is divisible by $p^s$.
Leibniz's rule and integral operator coefficients preserve the
joint congruences. The correction has increasing lowest degree,
and its coefficient at $Z^{p^{s-1}}$ is exactly $p$, since
$T_s^{(n)}(0)=1$. This proves convergence and the norm assertion.

For $s\geq2$ one has $U_s(Q)=U_{s-1}(Q^p)$, and hence
\begin{equation}\label{pt:twisted-frobenius-quotient}
 \frac{S_s^{(n)}(Q)}{S_{s-1}^{(n)}(Q^p)}
 =\frac{T_s^{(n)}(Q)}{T_{s-1}^{(n)}(Q^p)}.
\end{equation}
At the first level the quotients agree modulo $p$, because
$U_1\equiv1\pmod p$. The Dwork congruences for $T_s^{(n)}$
therefore give
\begin{equation}\label{pt:twisted-dwork}
 S_{s+1}^{(n)}(Q)S_{s-1}^{(n)}(Q^p)
 \equiv S_s^{(n)}(Q)S_s^{(n)}(Q^p)\pmod{p^s},\qquad s\geq1.
\end{equation}
Thus the joint equations, Dwork congruences, constant coefficient
one, and prescribed coefficientwise limit admit distinct
polynomial approximation sequences.

The comparison in Section~\ref{pt:section} concerns the normalized exact solution
$F_{A_s}^{(n)}$ at coupling $-A_s$, fixed by characteristic-zero
uniqueness. It computes its discrepancy from $T_s^{(n)}$,
including the sharp bound in terms of $B(\boldsymbol e)$ and the
leading residues.
\end{remark}

\section{Two- and Three-Particle Calculations}
\label{fl:section}
Let us consider the examples with $n=2$ and $n=3$ explicitly.

\subsection{The Two-Particle Case}
\label{ex:two-particle}
For $X_2=T^*\mathbb P^1$, write $Q=z_1/z_2$ and
$\theta=Q\partial_Q$. At the negative integer coupling $-m$,
the Gauss series of \eqref{cm:two-particle-vertex} terminates:
\begin{equation}\label{ex:two-particle-polynomial}
 F_m^{(2)}(Q)
 =\CT_u[(1-Q/u)^m(1-u)^m]
 =\sum_{d=0}^{m}\binom md^2Q^d.
\end{equation}
The constant term requires equal binomial indices, and their
recurrence gives the exact equation
\[
 [\theta^2-Q(\theta-m)^2]F_m^{(2)}=0.
\]
Thus $\widetilde T_s^{(2)}=T_s^{(2)}$ and
$\Delta_s^{(2)}=0$.

Under \eqref{ar:prime-assumptions}, put $P=p^s$ and
$m=\omega(P-1)$. Comparing the two quadratic operators gives
the complete residual at the arithmetic coupling:
\begin{equation}\label{ex:two-particle-residual}
 [\theta^2-Q(\theta+\omega)^2]T_s^{(2)}
 =-\omega P Q(2\theta+\omega-m)T_s^{(2)}
 \equiv0\pmod{p^s}.
\end{equation}
This identity involves no division by degree eigenvalues modulo
$p$. On ratio functions, the center integral is zero and the
second integral is $-H_{2,\omega}$, so this also gives the full
commuting congruence in two particles.

For $\omega=1/2$, $p=3$, the first levels are
\begin{equation}\label{ex:two-particle-levels}
 T_1^{(2)}=1+Q,\qquad
 T_2^{(2)}=1+16Q+36Q^2+16Q^3+Q^4.
\end{equation}
The first residual is $3Q(1-3Q)/4$, and
\[
 T_2^{(2)}\equiv(1+Q)(1+Q^3)
 =T_1^{(2)}(Q)T_1^{(2)}(Q^3)\pmod3.
\]
Consequently
$V_2(Q;1/2)\equiv\prod_{j\geq0}(1+Q^{3^j})\pmod3$,
the instance of the product formula already obtained in
\cite{SV23}. The stronger residual in the shifted
Schr\"odinger gauge is visible at the same level:
\begin{equation}\label{ex:two-particle-shifted}
 \phi_1=\frac{1+Q}{1-Q},\qquad
 H_{2,1/2}\phi_1=\frac94\frac{Q(1+Q)}{(1-Q)^3}
 \equiv0\pmod9.
\end{equation}
This is \eqref{pt:half-coupling-full} for two particles.

\subsection{The First Three-Particle Polynomial}
Write $u=x_{1,1}$, $v=x_{2,1}$, $w=x_{2,2}$ and put
\begin{equation}\label{fl:edge-polynomial}
 E=(u-v)(u-w)\prod_{j=1}^3(v-z_j)(w-z_j),\qquad
 G=z_2z_3^2,\qquad Q_1=z_1/z_2,\quad Q_2=z_2/z_3.
\end{equation}
In the chamber
$|z_1|<|v|<|u|<|z_2|<|w|<|z_3|$, factor
$E=uvw^3G\mathscr E$, where
\begin{equation}\label{fl:chamber-edges}
\begin{split}
 \mathscr E={}&(1-v/u)(1-u/w)(1-z_1/v)(1-v/z_2)\\
 &\cdot(1-v/z_3)(1-z_1/w)(1-z_2/w)(1-w/z_3).
\end{split}
\end{equation}
The period and its half-coupling approximants are
\begin{align}
 V_3(Q;\omega)&=\CT_{u,v,w}
 (1-v/w)^{2\omega}\mathscr E^{-\omega},\label{fl:period}\\
 T_s^{(3)}&=\CT_{u,v,w}(1-v/w)\mathscr E^m,
 \qquad m=(p^s-1)/2.\label{fl:finite-constant-term}
\end{align}
Equivalently,
\begin{equation}\label{fl:finite-extraction}
 T_s^{(3)}=G^{-m}[u^mv^mw^{3m+1}](w-v)E^m.
\end{equation}

For $p=3$, $s=1$, we have $m=1$. If
$F(t)=\prod_j(t-z_j)$, extraction in \eqref{fl:finite-extraction}
gives $[vw^4](v^2-w^2)F(v)F(w)=e_1(z)e_2(z)$.
Therefore
\begin{equation}\label{fl:p3-polynomial}
\begin{split}
 T_1^{(3)}&=(1+Q_1+Q_1Q_2)(1+Q_2+Q_1Q_2)\\
 &=1+Q_1+Q_2+3Q_1Q_2+Q_1^2Q_2+Q_1Q_2^2+Q_1^2Q_2^2.
\end{split}
\end{equation}
The product formula modulo $3$ becomes
\begin{equation}\label{fl:product}
 V_3(Q;1/2)\equiv\prod_{j\geq0}
 (1+Q_1^{3^j}+(Q_1Q_2)^{3^j})
 (1+Q_2^{3^j}+(Q_1Q_2)^{3^j})\pmod3.
\end{equation}

\subsection{The Energy Recurrence and Framing-Pole Cancellation}
Set $\theta_i=Q_i\partial_{Q_i}$ and write
$V_3=\sum_{d,e\geq0}c_{d,e}Q_1^dQ_2^e$, with $c_{0,0}=1$.
The quadratic equation in vertex gauge is equivalent to
\begin{equation}\label{fl:energy-recurrence}
\begin{split}
 (d^2-de+e^2)c_{d,e}={}&
 \omega\sum_{k=1}^{d}(2(d-k)-e+\omega)c_{d-k,e}\\
 &+\omega\sum_{k=1}^{e}(2(e-k)-d+\omega)c_{d,e-k}\\
 &+\omega\sum_{k=1}^{\min(d,e)}
 (d+e-2k+2\omega)c_{d-k,e-k}.
\end{split}
\end{equation}
In particular,
$c_{d,0}=c_{0,d}=((\omega)_d/d!)^2$ and
$c_{1,1}=2\omega^2(\omega^2-\omega+1)$.
At half coupling this gives
\begin{equation}\label{fl:initial-coefficients}
\begin{split}
 V_3(Q;1/2)={}&1+\frac14(Q_1+Q_2)
 +\frac9{64}(Q_1^2+Q_2^2)+\frac38Q_1Q_2\\
 &+\frac{25}{256}(Q_1^3+Q_2^3)
 +\frac{31}{256}(Q_1^2Q_2+Q_1Q_2^2)
 +O\bigl((Q_1,Q_2)^4\bigr).
\end{split}
\end{equation}
The recurrence is an exact characteristic-zero calculation.
The integral operator identities, rather than division by
$d^2-de+e^2$ modulo $p$, prove the finite-level congruences.

The localization formula gives the first mixed coefficient
directly and exhibits the cancellation of framing poles. Keep $(a_1,a_2,a_3)=(0,\delta,\xi)$ until the degree
sum is formed. In degree $(1,1)$ the two row-two splits give
\begin{align}
 A(\delta,\xi)&=
 \omega\frac{\omega+\delta}{1+\delta}
 \frac{\omega+\xi}{1+\xi}
 \frac{1+\delta-\omega}{\delta},\notag\\
 B(\delta,\xi)&=
 \omega^2\frac{\omega+\xi-\delta}{1+\xi-\delta}
 \frac{\omega+\delta-1}{\delta}.
 \label{fl:mixed-localization-terms}
\end{align}
Their residues at $\delta=0$ cancel. On the slice $\xi=0$,
the complete sum is
\begin{equation}\label{fl:mixed-localization-cancellation}
 A(\delta,0)+B(\delta,0)
 =2\omega^2\frac{\omega^2-\omega+1-\delta^2}{1-\delta^2}
 \longrightarrow2\omega^2(\omega^2-\omega+1).
\end{equation}
For $\omega=1/2$ this is $3/8$.
Proposition~\ref{dg:vertex-regularity} proves regularity of the
simultaneous framing specialization in every degree; the example
shows why the summands must be combined first.

\subsection{The Discrepancy at the First Arithmetic Levels}
Put $x=Q_1$, $y=Q_2$, and abbreviate
$C_{m,\gamma}=C_3((x,y);m,\gamma)$.
At half coupling define
\[
 B_m=C_{m,1},\qquad F_m=F_m^{(3)}=C_{m,-2m},\qquad
 D_m=F_m-B_m.
\]
Thus $D_m=\Delta_s^{(3)}$ when $m=(p^s-1)/2$.
The first mixed coefficient already distinguishes the families:
\begin{equation}\label{pt:first-difference}
 [xy]C_{m,\gamma}=2m^4+(2-\gamma)m^2,\qquad
 [xy]\Delta_s^{(3)}=p^sm^2.
\end{equation}
For $m=1$, $B_1$ is \eqref{fl:p3-polynomial} and
\begin{equation}\label{pt:first-polynomials}
 F_1=B_1+3xy,\qquad D_1=3xy,\qquad
 \mathcal L_{3,-1}D_1=D_1.
\end{equation}
The symmetric polynomial associated with this terminating vertex is
\[
 G_3(z)F_1=m_{(2,1,0)}(z)+6m_{(1,1,1)}(z)
 =P_{(2,1,0)}^{(-1)}(z).
\]
Here $m_\lambda$ is the sum of the distinct monomials obtained
by permuting the exponents of $z^\lambda$. This is the first
three-particle instance of \eqref{pt:jack-specialization}.
For $m=2$, corresponding to $p=5$, $s=1$, the full difference is
\begin{equation}\label{pt:second-polynomials}
 D_2=20xyB_1+30x^2y^2.
\end{equation}

At the next $3$-adic level, $P=9$ and $m=4$.
Use the weights of Lemma~\ref{pt:flow-weights}, with $r=v/w$.
Direct expansion of the finite edge product gives
\[
 W_1^{(4)}(3,3)=12160,\qquad
 W_2^{(4)}(3,3)=704,\qquad
 W_3^{(4)}(3,3)=16.
\]
Since $(1-r)^{-8}-(1-r)=9r+36r^2+120r^3+\cdots$,
\begin{equation}\label{pt:counterexample-coefficient}
 [x^3y^3]D_4
 =9\cdot12160+36\cdot704+120\cdot16
 =136704\equiv3\pmod9.
\end{equation}
The complete polynomial satisfies
\begin{equation}\label{pt:counterexample-polynomial}
 D_4\equiv3(xy)^3B_1\pmod9.
\end{equation}
For a direct derivation, only same-row indices $h\leq4$ can
contribute by Lemma~\ref{pt:flow-weights}. Reducing their
coefficients modulo $9$ gives
\[
 D_4/3\equiv\CT_{u,v,w}r^3(1-r)\mathscr E_3^4\pmod3.
\]
Put $\mathscr C=(1-r)\mathscr E_3$.
Lemma~\ref{ar:master-support}, applied at $p=3$, $s=1$,
says that the part of $\mathscr C$ whose internal exponents
are all divisible by $3$ is $B_1(x,y)(1-r^3)$.
Since $\mathscr E_3^4\equiv\mathscr E_3\sigma(\mathscr E_3)
\pmod3$, the preceding constant term becomes
\[
 \CT_{u,v,w}\bigl[\mathscr C\,\sigma(r\mathscr E_3)\bigr]
 =B_1(x,y)(xy)^3.
\]
Indeed, dividing the selected internal exponents by $3$ leaves
$\CT_{u,v,w}r(1-r)\mathscr E_3$, with external variables
replaced by $x^3,y^3$. Before this replacement the constant
term is $xy$: the index-one contribution is $xy$, and the
index-two contribution vanishes.
This proves \eqref{pt:counterexample-polynomial}.

In particular, $D_4$ is a nonzero joint kernel class
modulo $9$, while
$D_4/3\equiv(xy)^3B_1\pmod3$. The factor $(xy)^3$ is
annihilated by the Euler derivatives modulo $3$, as in
Remark~\ref{cg:uniqueness-scope}.

The coefficient $136704$ realizes the sharp bound of
Theorem~\ref{pt:sharpness} at the first higher $3$-adic level.
By \eqref{pt:sharp-face-coefficient}, the same coefficient occurs
on every consecutive three-particle subflag in higher rank.

\subsection{A Four-Particle Symmetric Example}
Consider the four-particle case with $p=3$, $s=1$ and
$\omega=1/2$, so $m=1$. The corresponding Jack polynomial
is indexed by the staircase partition $\delta=(3,2,1,0)$.
The two homogeneous polynomials have the monomial symmetric
expansions
\begin{align*}
 P_{3210}^{(-1)}
 &=m_{3210}+6m_{3111}+6m_{2220}+36m_{2211},\\
 \mathcal T_1
 &=m_{3210}+3m_{3111}+3m_{2220}+9m_{2211}.
\end{align*}
Here partitions are written without commas, and the $e_j$ below
are elementary symmetric polynomials in $z_1,\ldots,z_4$.
The master-polynomial extraction gives
\[
 \mathcal T_1=e_2(e_1e_3+e_4).
\]
Indeed, extraction in the first two internal rows gives
$e_1(w)e_2(w)$ in the three variables $w$ of row three.
Multiplication by $\Delta_3(w)$ gives the sum of alternants
with exponent sets $(0,2,4)$ and $(1,2,3)$. The last-row
extraction contributes $e_1e_2e_3$ and $e_2e_4$, respectively.
The Jack coefficients follow from the quadratic spectral
recursion on the four partitions dominated by $3210$.
Consequently
\begin{equation}\label{ex:rank-four-jack}
 P_{3210}^{(-1)}-\mathcal T_1
 =3m_{3111}+3m_{2220}+27m_{2211}.
\end{equation}
The partitions $3111$ and $2220$ are the two three-row
squeezes in \eqref{pt:three-row-squeeze}. The coefficients
of these terms have $3$-adic valuation one and attain the
sharp bound. The coefficient of $m_{2211}$ has valuation
three. This example also illustrates the comparison when
$p\leq n$, with no division by a symmetric-group orbit size.

\appendix

\section{Localization, the tRS Limit, and Zero Framing}
\label{dg:section}
We derive the scalar vertex equation from the K\"ahler-variable
eigenfunction theorem of \cite{K21} (Theorem~2.10), using the
normalization of \cite{Din21} (\S2.6, Definition~4 and
Propositions~1--2). The calculation fixes the degree variables,
the framing characters, and the scalar character accompanying the
bare vertex. A recursion over the local ring at zero framing then
proves cancellation of the framing poles in each complete degree
sum. Throughout this appendix, limits in the formal parameter
$\eta$ are coefficientwise in the K\"ahler variables.

\subsection{Localization and Degree Normalization}
\label{sec:geometry}
Choose the identity coordinate flag, set the source weight to
$\eps=1$, and write $a_1,\ldots,a_n$ for the additive framing
weights and $\omega=\hbar/\eps$ for the cotangent weight. For
$m\in\Z$, put
\[
 (b)_m=\begin{cases}
 \prod_{j=0}^{m-1}(b+j),&m\geq0,\\
 \displaystyle\prod_{j=1}^{-m}(b-j)^{-1},&m<0,
 \end{cases}
 \qquad R(b,c;m)=\frac{(b)_m}{(c)_m}.
\]
Exceptional arguments are substituted into the rational ratio,
not into the two Pochhammer symbols separately.

For nonnegative arrays $(d_{i,a})_{1\leq a\leq i<n}$, set
$d_{n,b}=0$ and $d_i=\sum_{a=1}^i d_{i,a}$. The localization
formula in this convention is
\begin{align}
 V_n(Q;a,\omega)
 ={}&\sum_{(d_{i,a})\geq0}\prod_{i=1}^{n-1}Q_i^{d_i}
 \prod_{i=1}^{n-1}\prod_{a,b=1}^{i}
 R(1-a_a+a_b,\omega-a_a+a_b;d_{i,a}-d_{i,b})\notag\\
 &\quad\cdot\prod_{i=1}^{n-1}\prod_{a=1}^{i}\prod_{b=1}^{i+1}
 R(\omega+a_b-a_a,1+a_b-a_a;d_{i,a}-d_{i+1,b}).
 \label{ge:localization}
\end{align}
Each degree coefficient is a finite sum of rational functions.
The same-index edge ratios enforce
$d_{i,a}\geq d_{i+1,a}$; this inequality must not be imposed on
arbitrary pairs of distinct indices. In particular, a negative
degree difference does not by itself make a localization term
vanish.

The specialization under consideration is
\begin{equation}\label{ge:vertex-limit}
 V_n(Q;\omega)=\left.V_n(Q;a,\omega)\right|_{a=0},
 \qquad V_n(0;\omega)=1.
\end{equation}
It is taken after summing all terms in each degree. Its existence
will follow from Proposition~\ref{dg:vertex-regularity}.

For the multiplicative formula, write $h$ for the cotangent
character and $\alpha_j$ for the framing characters, and set
\[
 (u;q)_m=\begin{cases}
 \prod_{j=0}^{m-1}(1-uq^j),&m\geq0,\\
 \displaystyle\prod_{j=1}^{-m}(1-uq^{-j})^{-1},&m<0,
 \end{cases}
 \qquad B_q(u;m)=\frac{(hu;q)_m}{(qu;q)_m}.
\]
The same rational-ratio convention applies at exceptional weights.
Let $z_i^{\rm raw}$ be the degree variables of the bare vertex in
\cite{KZ21} (Theorem~3.1), and put
\begin{equation}\label{cl:kahler-rescaling}
 Z_i=\frac{q^2}{h}z_i^{\rm raw},\qquad 1\leq i<n.
\end{equation}
Its identity component becomes
\begin{align}\label{cl:k-localization}
 \mathbf V_n^K(Z;\alpha,h,q)
 ={}&\sum_{(d_{i,a})\geq0}\prod_{i=1}^{n-1}Z_i^{d_i}
 \prod_{i=1}^{n-1}\prod_{a,b=1}^{i}
 B_q(\alpha_b/\alpha_a;d_{i,a}-d_{i,b})^{-1}\notag\\
 &\quad\cdot\prod_{i=1}^{n-1}\prod_{a=1}^{i}\prod_{b=1}^{i+1}
 B_q(\alpha_b/\alpha_a;d_{i,a}-d_{i+1,b}).
\end{align}
Arrays outside the degree chamber of the cited formula violate a
same-index inequality and have a zero edge ratio at generic
framing. Thus summing over all nonnegative arrays gives the same
series.

Here is the degree-dependent scalar in
\eqref{cl:kahler-rescaling}. The brace factor of the cited
formula satisfies
\[
 \{x\}_m=B_q(x^{-1};m)(-q^{1/2}h^{-1/2})^m.
\]
The inverse root factors have total brace exponent zero. The
edge and framing factors have total exponent
\[
 \sum_{i=1}^{n-1}\bigl((i+1)d_i-id_{i+1}\bigr)
 =2\sum_{i=1}^{n-1}d_i=N(d),\qquad d_n=0.
\]
The extra factor $q^{N(d)/2}$ in that formula therefore gives
\[
 q^{N(d)/2}(-q^{1/2}h^{-1/2})^{N(d)}
 =(q^2/h)^{\sum_i d_i},
\]
which proves the stated rescaling.

\subsection{The Cohomological Localization Limit}
\begin{lemma}\label{dg:pochhammer}
For $m\in\Z$ and generic $U,V,\eps$ with $\eps\ne0$, the ratio
below is regular at $\eta=0$, and
\begin{equation}\label{dg:pochhammer-limit}
 \left.
 \frac{(e^{\eta U};e^{\eta\eps})_m}
      {(e^{\eta V};e^{\eta\eps})_m}
 \right|_{\eta=0}
 =\frac{(U/\eps)_m}{(V/\eps)_m}.
\end{equation}
\end{lemma}

\begin{proof}
For $m\geq0$, cancel $\eta$ in each factor
$(1-e^{\eta(U+j\eps)})/(1-e^{\eta(V+j\eps)})$.
For $m=-b<0$, use instead
\[
 \prod_{j=1}^b
 \frac{1-e^{\eta(V-j\eps)}}{1-e^{\eta(U-j\eps)}}.
\]
Their constant terms are, respectively,
$(U+j\eps)/(V+j\eps)$ and
$(V-j\eps)/(U-j\eps)$, giving the assertion in both cases.
\end{proof}

\begin{proposition}\label{cl:generic-limit}
Set
\[
 q=e^{\eta\eps},\qquad h=e^{\eta H},\qquad
 \alpha_j=e^{\eta A_j},\qquad
 \omega=H/\eps,\qquad a_j=A_j/\eps,\qquad Z_i=Q_i.
\]
For generic additive weights, each coefficient of
\eqref{cl:k-localization} is regular at $\eta=0$, and
\begin{equation}\label{cl:vertex-limit}
 \left.\mathbf V_n^K(Q;e^{\eta A},e^{\eta H},e^{\eta\eps})
 \right|_{\eta=0}=V_n(Q;a,\omega).
\end{equation}
\end{proposition}

\begin{proof}
Lemma~\ref{dg:pochhammer} sends an edge ratio of degree $m$ to
$R(\omega+a_b-a_a,1+a_b-a_a;m)$ and an inverse root ratio to
its reciprocal, precisely as in \eqref{ge:localization}.
The rational products also evaluate the same-index zero factors,
which remain zero. There are finitely many arrays with fixed row
sums, so the limit passes through the complete degree sum.
\end{proof}

The eigenfunction normalization below uses $Z_i=(q/h)Q_i$.
Since $q/h\to1$, this change has the same cohomological limit.

\subsection{Centered tRS Operators and the Spectral Dictionary}
From now on put
$q=e^\eta$, $t=h=e^{\eta\omega}$,
$c=(n-1)/2$, and $\rho_i=(n+1-2i)/2$.
The commuting tRS operators are
\begin{equation}\label{dg:raw-family}
 \mathsf T_{r,\eta}=t^{r(r-1)/2}
 \sum_{|I|=r}\prod_{i\in I,\,j\notin I}
 \frac{tz_i-z_j}{z_i-z_j}\,e^{\eta D_I},
 \qquad D_I=\sum_{i\in I}D_i.
\end{equation}
Their centered form is
\begin{equation}\label{dg:normalized-trs}
 \mathsf A_{r,\eta}
 =t^{-rc}z^{-\omega\rho}\mathsf T_{r,\eta}z^{\omega\rho}
 =\sum_{|I|=r}a_I(\eta,Q)e^{\eta D_I},
\end{equation}
where
\begin{equation}\label{dg:chamber-family}
 a_I=
 \prod_{\substack{i<j\\i\in I,\ j\notin I}}
 \frac{1-tq_{ij}}{1-q_{ij}}
 \prod_{\substack{i<j\\i\notin I,\ j\in I}}
 \frac{1-t^{-1}q_{ij}}{1-q_{ij}}.
\end{equation}
Indeed, the chamber constant in the coefficient of
$\mathsf T_{r,\eta}$ is $t^{\sum_{i\in I}(i-1)}$;
conjugation and centering cancel it. Thus $a_I(\eta,0)=1$.

\begin{proposition}\label{dg:trigonometric-limit}
The first centered operator has expansion
\begin{equation}\label{dg:quadratic-expansion}
 \mathsf A_{1,\eta}
 =n+\eta\sum_iD_i+\eta^2\mathcal L_{n,\omega}+O(\eta^3).
\end{equation}
The first-order term vanishes on series in $Q$.
\end{proposition}

\begin{proof}
Write $u_{ij}=q_{ij}/(1-q_{ij})$ and
$a_{\{i\}}=1+\eta b_i+\eta^2c_i+O(\eta^3)$. Formula
\eqref{dg:chamber-family} gives
\[
 b_i=-\omega\sum_{j>i}u_{ij}
       +\omega\sum_{j<i}u_{ji},\qquad \sum_i b_i=0,
\]
and
\[
 \sum_i c_i=-\omega^2\sum_{i<j}u_{ij}
 +\omega^2\sum_{i<j<k}
 (u_{ij}u_{ik}-u_{ij}u_{jk}+u_{ik}u_{jk}).
\]
Since $q_{ik}=q_{ij}q_{jk}$, the expression in parentheses is
$-u_{ik}$. Hence
$\sum_i c_i=-\omega^2\sum_{i<j}(j-i)u_{ij}$.
Expanding the shifts gives \eqref{cg:normalized-operator}; also
$\sum_iD_i=0$ on ratio series.
\end{proof}

We now specify the input from the K\"ahler eigenfunction theorem.
The variables $u_i,\hbar$ of \cite{Din21} are our
$\alpha_i,h$, and its homogeneous variables $\zeta_i$ obey
$\zeta_i/\zeta_{i+1}=(h/q)Z_i$. Identifying $\zeta_i=z_i$
gives the parameter correspondence
\begin{equation}\label{dg:vertex-kahler-dictionary}
 Q_i=qz_i^{\rm raw},\qquad
 Z_i=\frac qh Q_i=\frac{q^2}{h}z_i^{\rm raw},\qquad
 \widehat V_n^K(Q;a)=
 \mathbf V_n^K\bigl((q/h)Q;\alpha,h,q\bigr),
 \quad \alpha_i=e^{\eta a_i}.
\end{equation}
Its constant coefficient is one.

Put
$\mathcal R_I=\prod_{i\in I,j\notin I}(hz_i-z_j)/(z_i-z_j)$.
In the notation of \cite{Din21} (Definition~4),
\[
 \mathsf D_r^{\rm Din}
 =h^{r(r+1)/2-rn}\sum_{|I|=r}\mathcal R_Ie^{\eta D_I}
 =h^{-r(n-1)}\mathsf T_{r,\eta}.
\]
The normalizing factor $\mathcal C_\eta$ in Proposition~1 of
that reference has shift ratios
\[
 \frac{p_i\mathcal C_\eta}{\mathcal C_\eta}
 =h^{n-i}\alpha_i^{-1},\qquad p_i=e^{\eta D_i},
\]
and its Proposition~2 gives
\[
 \mathsf D_r^{\rm Din}
 (\mathcal C_\eta\widehat V_n^K)
 =e_r(\alpha^{-1})\mathcal C_\eta\widehat V_n^K.
\]
Thus the spectral characters are the reciprocals of the framing
characters in localization. K\"ahler-independent factors of
$\mathcal C_\eta$ do not affect this equation.

For comparison with the raw family, one may use
\begin{equation}\label{dg:theta-character}
 \Gamma_\eta(z;a)=\prod_i
 \frac{\theta(h^{i-1}z_i;q)}{\theta(\alpha_i^{-1}z_i;q)},
 \qquad
 \frac{p_i\Gamma_\eta}{\Gamma_\eta}
 =\alpha_i^{-1}h^{1-i},
\end{equation}
where $\theta(qx;q)=-x^{-1}\theta(x;q)$.
These shift ratios are also realized by the formal character
\begin{equation}\label{dg:formal-vertex-character}
 \Gamma_\eta^{\rm form}
 =\prod_i z_i^{-a_i-\omega(i-1)}
 =\left(\prod_i z_i\right)^{-\omega c}z^{\omega\rho-a}.
\end{equation}
Conjugation depends only on shift ratios, so no asymptotics of a
theta function are needed.

\begin{proposition}\label{dg:kahler-input}
For $1\leq r\leq n$, the normalized vertex satisfies
\begin{equation}\label{dg:vertex-monodromy}
 \mathsf A_{r,\eta}(z^{-a}\widehat V_n^K)
 =e_r(\alpha_1^{-1},\ldots,\alpha_n^{-1})
 z^{-a}\widehat V_n^K.
\end{equation}
Equivalently, with $\alpha_I=\prod_{i\in I}\alpha_i$,
\begin{equation}\label{dg:weighted-vertex-equation}
 \sum_{|I|=r}\alpha_I^{-1}a_I(\eta,Q)
 e^{\eta D_I}\widehat V_n^K
 =e_r(\alpha^{-1})\widehat V_n^K.
\end{equation}
\end{proposition}

\begin{proof}
For $|I|=r$, conjugation of $\mathsf D_r^{\rm Din}$ by
$\mathcal C_\eta$ gives the coefficient
\[
 h^{r(r+1)/2-rn}\alpha_I^{-1}
 h^{\sum_{i\in I}(n-i)}\mathcal R_I
 =\alpha_I^{-1}
 h^{r(r-1)/2+\sum_{i\in I}(1-i)}\mathcal R_I
 =\alpha_I^{-1}a_I.
\]
The last equality follows by factoring each cross pair according
to its order. Consequently
\begin{equation}\label{dg:theta-conjugation}
 \mathcal C_\eta^{-1}\mathsf D_r^{\rm Din}\mathcal C_\eta
 =\Gamma_\eta^{-1}\mathsf T_{r,\eta}\Gamma_\eta
 =\sum_{|I|=r}\alpha_I^{-1}a_Ie^{\eta D_I}.
\end{equation}
This proves \eqref{dg:weighted-vertex-equation} from the cited
eigenfunction theorem. A shift in $r$ variables multiplies the
center factor of \eqref{dg:formal-vertex-character} by $h^{-rc}$,
so
$\Gamma_\eta^{-1}\mathsf T_{r,\eta}\Gamma_\eta
=z^a\mathsf A_{r,\eta}z^{-a}$, proving
\eqref{dg:vertex-monodromy}.
\end{proof}

For $n=2$, the normalized $K$-theoretic vertex is
\[
 \widehat V_2^K(Q;a)={}_2\phi_1
 \left(h,h\alpha_2/\alpha_1;
 q\alpha_2/\alpha_1;q;(q/h)Q\right).
\]
Its coefficient ratio gives the two reciprocal framing weights
in \eqref{dg:weighted-vertex-equation}.

\subsection{The Local-Ring Argument at Zero Framing}
On series in $Q$ define
\begin{equation}\label{dg:generic-divided-operator}
 \mathsf N_\eta(a)=\eta^{-2}
 \left[\sum_i e^{-\eta a_i}a_{\{i\}}(\eta,Q)e^{\eta D_i}
 -\sum_i e^{-\eta a_i}\right].
\end{equation}
The order-zero term in brackets vanishes. The order-one term is
$\sum_iD_i$, hence also vanishes on ratio series; the divided
operator is therefore regular at $\eta=0$. The expansion above
gives
\begin{equation}\label{dg:generic-vertex-equation}
 \left.\mathsf N_\eta(a)\right|_{\eta=0}
 =\mathcal L_{n,\omega,a}
 :=\mathcal L_{n,\omega}-\sum_i a_iD_i
 +\omega\sum_{i<j}\frac{q_{ij}}{1-q_{ij}}(a_i-a_j).
\end{equation}
Propositions~\ref{dg:kahler-input} and~\ref{cl:generic-limit}
then imply, for generic framing,
\begin{align}\label{dg:generic-cohomological-vertex}
 \mathcal L_{n,\omega,a}V_n(Q;a,\omega)&=0,\\
 \mathcal L_{n,\omega}(z^{-a}V_n(Q;a,\omega))
 &=\frac12\sum_i a_i^2\,z^{-a}V_n(Q;a,\omega).\notag
\end{align}
The corresponding Schr\"odinger function is
$z^{-a}J_nV_n(Q;a,\omega)$.

\begin{proposition}\label{dg:vertex-regularity}
Every coefficient of $V_n(Q;a,\omega)$ is regular at $a=0$.
More precisely, each coefficient of $\widehat V_n^K(Q;a)$ lies in
$R[[\eta]]$, where
\[
 R=\Q(\omega)[a_1,\ldots,a_n]_{(a_1,\ldots,a_n)}.
\]
The cohomological and zero-framing specializations commute:
\begin{equation}\label{dg:joint-regularity}
 \left.\widehat V_n^K(Q;a)\right|_{\eta=0,\ a=0}
 =\left.V_n(Q;a,\omega)\right|_{a=0}
 =\left.\widehat V_n^K(Q;0)\right|_{\eta=0}.
\end{equation}
Moreover $[Q^{\boldsymbol d}]V_n(Q;\omega)\in\Q[\omega]$.
\end{proposition}

\begin{proof}
Construct a normalized solution of $\mathsf N_\eta(a)W=0$
over $R[[\eta]]$. Put $d_0=d_n=0$ and
$\nu_i=d_i-d_{i-1}$. Since $a_I(\eta,0)=1$, the diagonal
coefficient in degree $\boldsymbol d$ is
\begin{equation}\label{dg:joint-recursion-diagonal}
 E_{\boldsymbol d}(a,\eta)
 =\frac{\sum_i e^{-\eta a_i}(e^{\eta\nu_i}-1)}{\eta^2}
 =\frac12\sum_i\nu_i^2-\sum_i a_i\nu_i+O(\eta).
\end{equation}
Here $\sum_i\nu_i=0$. For $\boldsymbol d\ne0$, its constant
term at $\eta=a=0$ is the positive integer
$\mathcal E(\boldsymbol d)$ of \eqref{cg:indicial-form}.
Thus $E_{\boldsymbol d}$ is a unit in $R[[\eta]]$.
All off-diagonal terms involve smaller total K\"ahler degree.
Starting with $w_{\boldsymbol0}=1$, division by these units
therefore constructs a unique
$W=\sum_{\boldsymbol d}w_{\boldsymbol d}Q^{\boldsymbol d}$ with
$w_{\boldsymbol d}\in R[[\eta]]$. The degree-zero equation is
identically zero.

Over $K=\Q(\omega,a_1,\ldots,a_n)$, the same recursion gives
uniqueness in $K((\eta))[[Q]]$.
Proposition~\ref{dg:kahler-input} places the generic localization
vertex in this equation, and Proposition~\ref{cl:generic-limit}
gives its coefficientwise expansion at $\eta=0$.
It must therefore equal $W$. In particular, the complete degree
sums, although their summands may have framing poles, belong to
$R[[\eta]]$.

The cohomological limit identifies $W|_{\eta=0}$ with
$V_n(Q;a,\omega)$. Evaluation at $\eta=0$ and evaluation at
$a=0$ commute on $R[[\eta]]$, proving
\eqref{dg:joint-regularity}.

Finally let
$R_0=\Q[a_1,\ldots,a_n]_{(a_1,\ldots,a_n)}$.
The coefficients of $\mathsf N_\eta(a)$ are polynomial in
$\omega$ over $R_0$, whereas the diagonal
$E_{\boldsymbol d}(a,\eta)$ is independent of $\omega$ and has
its inverse in $R_0[[\eta]]$.
The recursion consequently gives
$w_{\boldsymbol d}\in(R_0[\omega])[[\eta]]$.
Evaluating at $\eta=a=0$ proves the last assertion and permits
every characteristic-zero specialization of the coupling.
\end{proof}

\begin{proposition}\label{dg:vertex-equation}
The normalized scalar vertex \eqref{ge:vertex-limit} exists in
$\Q[\omega][[Q]]$ and satisfies
\begin{equation}\label{dg:vertex-period-identification}
 V_n(Q;\omega)=F_{n,\omega}(Q).
\end{equation}
In particular, $\mathcal L_{n,\omega}V_n=0$ and
$H_{n,\omega}(J_nV_n)=0$.
\end{proposition}

\begin{proof}
Set $a=0$ in \eqref{dg:generic-cohomological-vertex} using
Proposition~\ref{dg:vertex-regularity}, and apply the normalized
formal uniqueness of Proposition~\ref{cg:formal-uniqueness}.
Conjugation by $J_n$ gives the Schr\"odinger equation.
\end{proof}

\section{Adjacent-Rank Kernels and Chamber Constant Terms}
\label{tc:section}
\label{tc:kernel-appendix}
This appendix proves the quadratic equation for the master
period and identifies the normalization of the commuting
matching operators. The matching formula also identifies the physical
coupling dependence used in \eqref{pt:half-coupling-full}.

\subsection{The Kernel and Its Chamber Expansion}
For a row of $r$ variables put
$\Delta_r(v)=\prod_{a<b}(v_b-v_a)$, with $\Delta_1=1$.
For adjacent rows $z=(z_1,\ldots,z_r)$ and
$y=(y_1,\ldots,y_{r-1})$, define, up to a constant branch factor,
\begin{equation}\label{tc:kernel-definition}
 \mathscr K_r(z,y)
 =\Delta_r(z)^\kappa\Delta_{r-1}(y)^\kappa
   \prod_{a=1}^r\prod_{b=1}^{r-1}(z_a-y_b)^{-\kappa}
   \prod_{b=1}^{r-1}y_b^\kappa.
\end{equation}
We always use the interlacing chamber
$|z_1|<|y_1|<|z_2|<\cdots<|y_{r-1}|<|z_r|$.
Factoring the differences cancels every monomial power.
We choose the constant phase so that the remaining product of
binomial series has leading coefficient one.

\begin{lemma}\label{tc:chamber-ring}
Choose a total ordering of all internal and external radii
compatible with \eqref{ar:interlacing}, and write the variables
as $v_1,\ldots,v_L$, where $L=n(n+1)/2$, $v_1=z_1$ and
$v_L=z_n$. Put $w_j=v_j/v_{j+1}$.
The normalized kernels and the prescribed expansions of their
logarithmic derivatives belong to
$\mathscr A_K=K[[w_1,\ldots,w_{L-1}]]$.
Internal constant-term extraction defines a map
$\CT_x:\mathscr A_K\to K[[Q]]$ with
\[
 \CT_x D_{r,a}f=0,\qquad
 D_i\CT_x f=\CT_xD_if.
\]
It commutes with multiplication by external chamber series.
For $\kappa\in\Z_p$ the kernels are integral units; the same
assertions hold over $\Z_p$.
\end{lemma}

\begin{proof}
Every smaller-to-larger ratio is a nonconstant monomial in the
$w_j$. The normalized kernels are products of binomial series
in these monomials, and their logarithmic derivatives are rational
series of the same kind. For a monomial $w^{\boldsymbol e}$,
zero exponent of an internal variable says that its two adjacent
$w$ exponents agree. Thus $[Q^{\boldsymbol d}]\CT_x$ selects
the unique exponent vector that equals $d_j$ throughout the gap
from $z_j$ to $z_{j+1}$. This proves well-definedness and the
derivative identities. External multiplication preserves
internal exponents. Integrality follows from
$\binom\alpha j\in\Z_p$ for $\alpha\in\Z_p$.
\end{proof}

\subsection{The Quadratic Kernel Identity}
\begin{lemma}\label{tc:kernel}
For $r\geq2$,
\begin{equation}\label{tc:kernel-equation}
 (H_{r,\kappa}^{z}-H_{r-1,\kappa}^{y})\mathscr K_r(z,y)=0.
\end{equation}
After division by the kernel, this is a rational identity
polynomial in $\kappa$.
\end{lemma}

\begin{proof}
This is the adjacent-rank trigonometric specialization of the
kernel identity in \cite{Lan06}; we give the calculation in
the normalization used here.
Combine both rows into variables $v_a$, with $\sigma_a=1$ on
the $r$ upper variables and $\sigma_a=-1$ on the $r-1$ lower
variables. Up to a constant,
\[
 \mathscr K_r=
 \prod_{a<b}\left(\frac{v_a-v_b}{\sqrt{v_av_b}}\right)
                   ^{\kappa\sigma_a\sigma_b}.
\]
The induced monomial powers are zero on upper variables and
$\kappa$ on lower variables, as in
\eqref{tc:kernel-definition}. Put
\[
 R_{ab}=\frac{v_a+v_b}{v_a-v_b},\qquad
 V_{ab}=\frac{v_av_b}{(v_a-v_b)^2},\qquad
 \ell_a=\kappa^{-1}D_a\log\mathscr K_r
       =\frac12\sum_{b\ne a}\sigma_a\sigma_bR_{ab}.
\]
Then $D_a\ell_a=-\sum_{b\ne a}\sigma_a\sigma_bV_{ab}$.
The identities
\[
 R_{ab}^2=1+4V_{ab},\qquad
 R_{ab}R_{ac}-R_{ab}R_{bc}+R_{ac}R_{bc}=1
\]
give
\begin{align*}
 \frac12\sum_a\sigma_aD_a\ell_a
 &=-\sum_{\substack{a<b\\\sigma_a=\sigma_b}}\sigma_aV_{ab},\\
 \frac12\sum_a\sigma_a\ell_a^2
 &=\sum_{\substack{a<b\\\sigma_a=\sigma_b}}\sigma_aV_{ab}.
\end{align*}
In the second formula the constant pair contribution is
$(r-1)/4$ and the triple contribution is $-(r-1)/4$;
the latter is the coefficient of $X^3$ in
$(1+X)^r(1-X)^{r-1}$, divided by four.
The signed kinetic term is therefore
$\kappa(\kappa-1)\sum_{\sigma_a=\sigma_b}\sigma_aV_{ab}$,
which cancels the signed potential. The divided identity is
polynomial in $\kappa$, so the calculation also covers
$\kappa=0$.
\end{proof}

Write $x^{(r)}=(x_{r,1},\ldots,x_{r,r})$, with $x^{(n)}=z$,
and put $D_{r,a}=x_{r,a}\partial_{x_{r,a}}$,
$\mathscr K_1=1$ and
\[
 \mathscr F_r=\prod_{j=1}^{r-1}
       \mathscr K_{j+1}(x^{(j+1)},x^{(j)}),\qquad
 \mathscr F_1=1.
\]

\begin{proposition}\label{tc:divergence}
The product of kernels satisfies
\begin{equation}\label{tc:divergence-equation}
 H_{n,\kappa}^{z}\mathscr F_n
 =\sum_{r=1}^{n-1}\sum_{a=1}^r
       D_{r,a}(\mathscr F_n b_{r,a}),
\end{equation}
where
\begin{equation}\label{tc:divergence-coefficients}
 b_{r,a}=\frac\kappa2\left(
 1-\sum_{b=1}^{r+1}\frac{x_{r,a}}{x_{r,a}-x_{r+1,b}}
   +\sum_{b=1}^{r-1}\frac{x_{r,a}}{x_{r,a}-x_{r-1,b}}
 \right).
\end{equation}
The last sum is absent for $r=1$.
\end{proposition}

\begin{proof}
For the kinetic normalization used here,
\[
 fHK-KHf=\frac12\sum_aD_a(fD_aK-KD_af).
\]
Apply this to $f=\mathscr F_{r-1}$ and $K=\mathscr K_r$,
and use Lemma~\ref{tc:kernel}. The lower-row divergences
pass through the next kernel. The new primitive has coefficient
$\frac12D_{r-1,a}\log(\mathscr K_r/\mathscr K_{r-1})$.
The same-row Vandermonde terms cancel in this logarithmic
derivative, leaving \eqref{tc:divergence-coefficients}.
Induction from $H_{1,\kappa}1=0$ proves the formula.
\end{proof}

Multiplying the normalized kernels gives
\begin{equation}\label{tc:full-density}
 \mathscr F_n
 =J_{n,\kappa}
  \prod_{r=1}^{n-1}\prod_{a<b}
       (1-x_{r,a}/x_{r,b})^{2\kappa}
  \prod_{r=1}^{n-1}\prod_{a=1}^r\prod_{b=1}^{r+1}
       (1-\rho_{r,a,b})^{-\kappa}.
\end{equation}
Indeed, before chamber factorization the product is, up to a
constant phase,
\[
 \Delta_n(z)^\kappa M^\kappa\Delta^{2\kappa}E^{-\kappa}.
\]
The monomial factors cancel as in the normalized master density.
In particular,
\begin{equation}\label{tc:period-constant-term}
 \CT_x\mathscr F_n=J_{n,\kappa}\mathcal V_n.
\end{equation}

\begin{proposition}\label{tc:period-equation}\label{cg:period-equation}
The master period satisfies
\[
 H_{n,\kappa}(J_{n,\kappa}\mathcal V_n)=0,
 \qquad \mathcal L_{n,\kappa}\mathcal V_n=0.
\]
Hence $\mathcal V_n=F_{n,\kappa}$.
\end{proposition}

\begin{proof}
Take the internal constant term in
\eqref{tc:divergence-equation}. Its right side vanishes by
Lemma~\ref{tc:chamber-ring}; its left side becomes the
Hamiltonian applied to \eqref{tc:period-constant-term}.
Conjugation gives the vertex-gauge equation. The period has
constant coefficient one, so
Proposition~\ref{cg:formal-uniqueness} identifies it.
\end{proof}

\subsection{Matching Operators and Spectral Identification}
Let $\mathsf M_{r,k}(\kappa)$ denote the right side of the
matching formula \eqref{tc:matching-formula}, with $n$ replaced
by $r$. Set $\mathsf M_{r,0}=1$ and
$\mathsf M_{r,k}=0$ for $k>r$.

\begin{proposition}\label{tc:matching}
The operators $\mathsf M_{r,k}$ commute, have chamber constant
$e_k(D)$, and have integral chamber coefficients for
$\kappa\in\Z_p$. Their first two members are
\[
 \mathsf M_{r,1}=\sum_iD_i,\qquad
 \mathsf M_{r,2}=\frac12\left(\sum_iD_i\right)^2-H_{r,\kappa}.
\]
\end{proposition}

\begin{proof}
Use the commuting operators in \cite{HR15} (equations
(1.10)--(1.12) and (2.7)--(2.8)). Set $z_i=e^{2t_i}$ and
$\lambda=\kappa$ there. The derivative becomes $-2iD_i$,
and the pair factor becomes $-4g_\kappa v_{ij}$.
Dividing the order-$k$ operator by $(-2i)^k$ gives the formula
used here. A term with $\ell$ pairs occurs
$(r-k)!2^\ell\ell!(k-2\ell)!$ times in their permutation
sum, canceling the factorial prefactor before reduction modulo
$p$. Their commutation identities are rational identities
polynomial in $\kappa$, so hold at every coupling.
Finally $v_{ij}=q_{ij}/(1-q_{ij})^2$ has positive chamber
degree and integral coefficients. This proves the chamber
constant, the first two formulas and integrality.
\end{proof}

The matching family can be identified with the determinant
family using generic spectral parameters. This avoids a second
derivation of the joint period equations by higher kernel identities.

\begin{proposition}\label{tc:matching-identification}
For every $r$ and $0\leq k\leq r$,
$\mathsf M_{r,k}(\kappa)=\mathsf I_{r,k}(\kappa)$.
\end{proposition}

\begin{proof}
For independent framing parameters $a_i$, the normalized quadratic
solution $\psi_a=z^{-a}(1+O(Q))$ is unique: its positive-degree
diagonals $\mathcal E(\boldsymbol e)-a\cdot\nu(\boldsymbol e)$
are nonzero over $\Q(\kappa,a)$.
Proposition~\ref{tc:matching} and commutation with the Hamiltonian
therefore give $\mathsf M_{r,k}\psi_a=e_k(-a)\psi_a$.
Proposition~\ref{dg:full-family} gives the same action for
$\mathsf I_{r,k}$.

Write their difference in normal order as
$\sum_{\boldsymbol e\geq0}Q^{\boldsymbol e}C_{\boldsymbol e}(D)$.
Its chamber constant is zero. A nonzero coefficient of least
total degree would give $C_{\boldsymbol e}(-a)=0$ on applying
the difference to $\psi_a$, contrary to independence of the
$a_i$. Thus the families agree over $\Q(\kappa)$.
Polynomial dependence of each chamber coefficient on $\kappa$
extends the identity to every characteristic-zero specialization.
\end{proof}

\end{document}